\documentclass[english]{article}

\usepackage[T1]{fontenc}
\usepackage[latin1]{inputenc}
\usepackage{amsmath}
\usepackage{graphicx}

\newcommand{\miolabel}[1]{\label{#1}\fbox{#1}}

\renewcommand{\miolabel}{\label}

\usepackage{amsmath,amssymb}
\usepackage{babel}
\begin{document}

\title{3D simulations of early blood vessel formation}

\author{Fausto Cavalli\footnotemark[2] 
  \and Andrea Gamba \footnotemark[3]
  \and Giovanni Naldi \footnotemark[2]
  \and Matteo Semplice\footnotemark[2]
  \and Donatella Valdembri\footnotemark[4]
  \and Guido Serini\footnotemark[4]
}

\renewcommand{\thefootnote}{\fnsymbol{footnote}}
\footnotetext[2]{Dipartimento di Matematica, Universit\`a di Milano,
Via Saldini 50, I-20100 Milano, Italy.}
\footnotetext[3]{Dipartimento di Matematica, Politecnico di Torino, 
 Corso Duca degli Abruzzi 24, 10129 Torino, Italy.}
\footnotetext[4]{Department of Oncological Sciences and Division of
  Molecular Angiogenesis, 
  Institute for Cancer Research and Treatment, University of Torino
  School of Medicine, Strada Provinciale 142 Km 3.95, 10060 Candiolo
  (TO), Italia.}
\renewcommand{\thefootnote}{\arabic{footnote}}

\maketitle

\begin{abstract}
Blood vessel networks form by spontaneous aggregation of individual
cells migrating toward vascularization sites (vasculogenesis).  
A successful theoretical model of two dimensional experimental vasculogenesis 
has been recently proposed, showing the
relevance of percolation concepts and of cell cross-talk (chemotactic
autocrine loop) to the understanding of this self-aggregation
process. Here we study the natural 3D extension of the 
computational  model proposed earlier, which is relevant for the investigation of
the genuinely threedimensional process of vasculogenesis in vertebrate
embryos.  The computational model is based on a multidimensional Burgers
equation coupled with a reaction diffusion equation for a chemotactic
factor and a mass conservation law. 
The numerical approximation of the computational model is obtained by
high order relaxed schemes.  
Space and time  discretization are performed by using TVD
schemes and, respectively, IMEX schemes.  
Due to the computational costs of realistic simulations, we have
implemented the numerical algorithm on a cluster for parallel
computation. 
Starting from initial conditions mimicking the experimentally observed
ones, numerical simulations produce network-like structures qualitatively
similar to those observed in the early stages of \emph{in vivo} vasculogenesis.
We develop the computation of critical percolative indices as a robust measure
of the network geometry as a first step towards the comparison of computational and experimental data. 
\end{abstract}

\section{Introduction}
In recent years, biologists have collected many qualitative and
quantitative data on the behavior of microscopic components of living
beings. We are, however, still far from understanding in detail how
these microscopic components interact to build functions which are
essential for life. A problem of particular interest which has been
extensively investigated is the formation of patterns in biological
tissues \cite{KM94}.  Such patterns often show self-similarity and
scaling laws \cite{Man88} similar to those emerging in the physics of
phase transitions \cite{Sta87}.

The vascular network \cite{WBE97,WBE99} is a typical example of
natural structure characterized by non trivial scaling laws. In recent
years many experimental investigations have been performed on the
mechanism of blood vessel formation \cite{Car00} both in living beings
and in \textit{in vitro} experiments.  Vascular
networks form by spontaneous aggregation of individual cells
travelling toward vascularization sites (vasculogenesis). A successful
theoretical model of two dimensional experimental vasculogenesis 
has been recently proposed, showing the
relevance of percolation concepts and of cell cross-talk (chemotactic
autocrine loop) to the understanding of this self-aggregation
process.

Theoretical and computational modelling is useful in testing
biological hypotheses in order to explain which kind of coordinated
dynamics gives origin to the observed highly structured tissue patterns.
One can develop computational models based on simple
dynamical principles and test whether they are able to reproduce the
experimentally observed features. If the basic dynamical principles
are correctly chosen, computational experiments allow to observe the
emergence of complex structures from a multiplicity of interactions
following simple rules.

Apart from the purely theoretical interest, reproducing biological
dynamics by computational models allows to identify those biochemical
and biophysical parameters which are the most important in driving the process. This
way, computational models can produce a deeper understanding of biological
mechanisms, which in principle may end up having relevant practical
consequences. 
It is worth noticing here that a complete understanding of the
vascularization process is possible only if it is considered in its
natural threedimensional setting (\cite{Abb03,CPS+01}).

In this paper we illustrate computational results regarding the simulation
of vascular network formation
in a threedimensional environment.
We consider the threedimensional version of the model
proposed in \cite{GAC+03,SAG+03}. The model is based on a
Burgers-like equation, a well studied paradigm in the theory of
pattern formation, integrated with a feedback term describing the
chemotactic autocrine loop.  The numerical evolution of the
computational model starting from initial conditions mimicking the
experimentally observed ones produces network-like structures
qualitatively similar to those observed in the early stages of
\emph{in vivo} vasculogenesis.

Since in the long run we are interested in developing quantitative
comparison between experimental data and theoretical model, we start
by selecting a set of observable quantities providing robust
quantitative information on the network geometry. The lesson learned
from the study of twodimensional vasculogenesis is that percolative
exponents \cite{SA94} are an interesting set of such observables, so
we test the computation of percolative exponents on simulated network
structures.

A thorough quantitative comparison of the geometrical properties of experimental
and computational network structures will become possible as soon as
an adequate amount of experimental data, allowing proper statistical
computation, will become available.

The paper is organized as follows. Section 2 summarizes some
background knowledge on the biological problem of vascular network
formation. 
Section 3 is a short review of the properties of the
model introduced in \cite{GAC+03,SAG+03}. In Section 4 the numerical
approximation technique for the model is described. In Section 5 we
describe the qualitative properties of simulated network
structures and  present the results of the computation of
the exponents of the percolative transition. Finally, in the
Conclusions, we point out at predictable developments of our research.

\section{Biological background}

To supply tissues with nutrients in an optimal way, vertebrates have
developed a hierarchical vascular system which terminates in a network
of size-invariant units, \textit{i.e.} capillaries. Capillary networks
characterized by intercapillary distances ranging from $50$ to $300\,\mu\mathrm{m}$
are essential for optimal metabolic exchange~\cite{GH00}.

Capillaries are made of endothelial cells. Their growth is essentially
driven by two processes: vasculogenesis and angiogenesis~\cite{Car00}.
Vasculogenesis consists of local differentiation of precursor
cells to endothelial ones, that assemble into a vascular network by
directed migration and cohesion. Angiogenesis is essentially characterized
by sprouting of novel structures and their remodelling.

In twodimensional assays, the process of formation of a vascular
network starting from randomly seeded cells can be accurately tracked
by videomicroscopy \cite{GAC+03} and it is observed to proceed along
three main stages: \textit{i)} migration and early network formation,
\textit{ii)} network remodelling and \textit{iii)} differentiation in
tubular structures. During the first phase, which is the most
important for determining the final geometrical properties of the
structures, cells migrate over distances which are an order of
magnitude larger than their radius and aggregate when they adhere with
one of their neighbours.  An accurate statistics of individual cells
trajectories has been presented in \cite{GAC+03}, showing that, in the
first stage of the dynamics, cell motion has marked directional
persistence, pointing toward zones of higher cell
concentration. This indicates that cells communicate 
through the emission of soluble chemical factors that diffuse (and
degrade) in the surrounding medium, moving toward
the gradients of this chemical field.  Cells behave like not-directly
interacting particles, the interaction being mediated by the release
of soluble chemotactic factors. Their dynamics is well reproduced by the
theoretical model proposed in \cite{GAC+03}.

The lessons learned from the study of in vitro vasculogenesis is thus
that the formation of experimentally observed structures can be explained as the
consequence of cell motility and of cell cross-talk mediated by the
exchange of soluble chemical factors (chemotactic autocrine loop).
The theoretical model also shows that the main factors determining the qualitative
properties of the observed vascular structures are the available cell
density and the diffusivity and half-life of the soluble chemical
exchanged. It seems that only the dynamical rules followed by the
individual cell are actually encoded in the genes. 
The interplay of these simple dynamical 
rules with the geometrical and physical properties of the environment produces the highly
structured final result. 

At the moment, no direct observation of the chemotactic autocrine loop regulating vascular 
network formation is available, although
several indirect biochemical observations point to it, so, the
main evidence in this sense still comes from the theoretical
analysis of computational models. 

Several major developments in threedimensional cell culture and in
cell and tissue imaging allow today to observe in real time the
mechanisms of cell migration and aggregation in threedimensional
settings \cite{Fri04b,KPS06}.

In the embryo, endothelial cells are produced and migrate in a
threedimensional scaffold, the extracellular matrix. Migration is
actually performed through a series of biochemical processes, such as
sensing of chemotactic gradients, and of mechanical operations, such
as extensions, contractions, and degrading of the extracellular matrix
along the way.

The evidence provided by twodimensional experimental
vasculogenesis suggests that cell motion can be directed by an autocrine loop of
soluble chemoattractant factors 
also in the real threedimensional environment.

As a sample of typical vascular structures that are observe in a
threedimensional setting in the early stages of development of a
living being, we include here $(750\,\mu\mathrm{m)^{2}}$ images of
chick embryo brain at different development stages
(Fig.~\ref{fig:vascular}).
At an early stage (about 52-64 hours) one observes a typical immature
vascular network formed by vasculogenesis and characterized by a high
density of similar blood vessels (Fig.~\ref{fig:vascular}A). At the
next stage (70-72 hours) we observe initial remodelling of the vascular
network (Figs.~\ref{fig:vascular}B,C).  Remodeling becomes more
evident when the embryo is 5 days old, when blood vessels are
organized in a mature, hierarchically organized vascular tree
(Fig.~\ref{fig:vascular}D).

\begin{figure}
\textsf{\textbf{A}} \includegraphics[%
  scale=0.3]{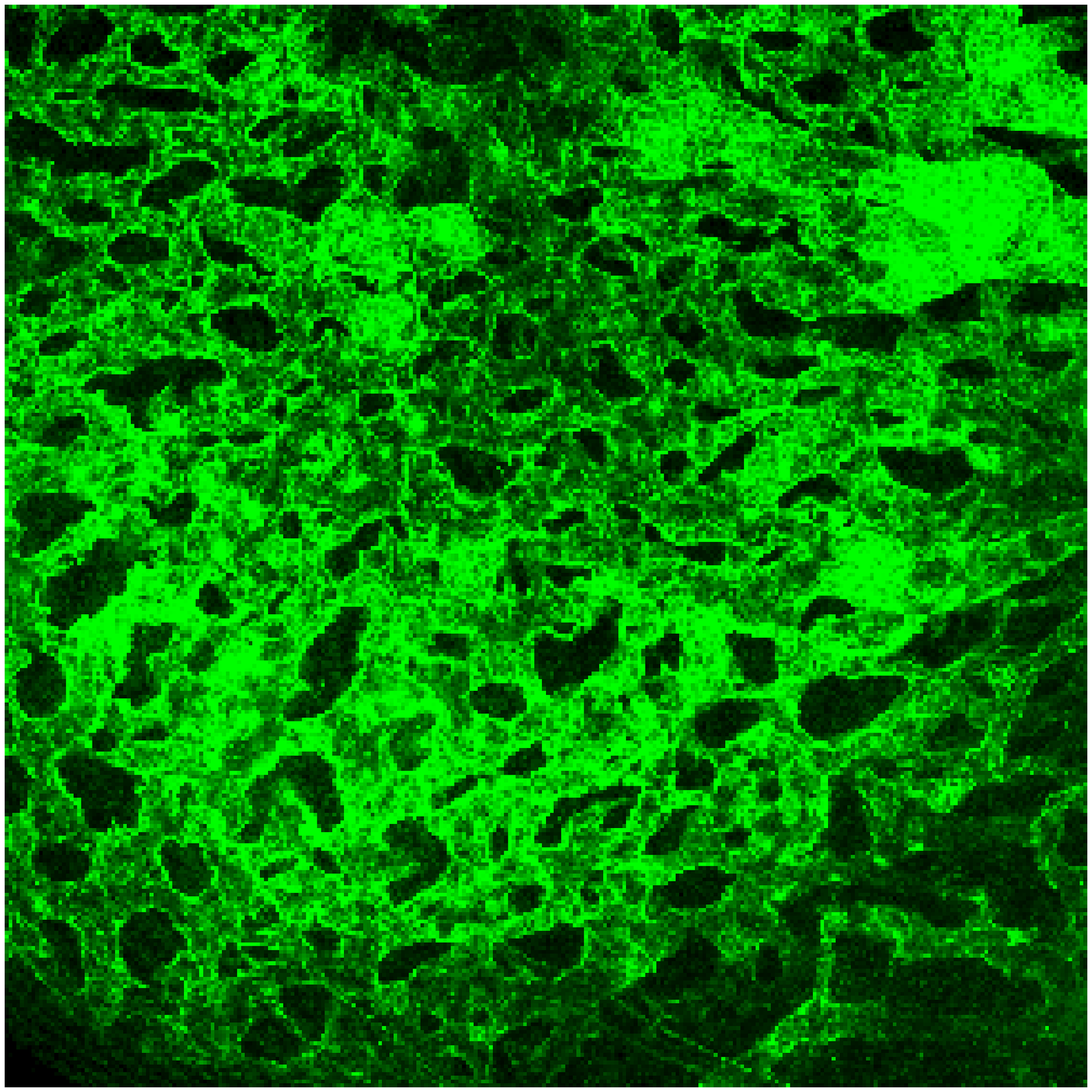} \textsf{\textbf{B}}\includegraphics[%
  scale=0.3]{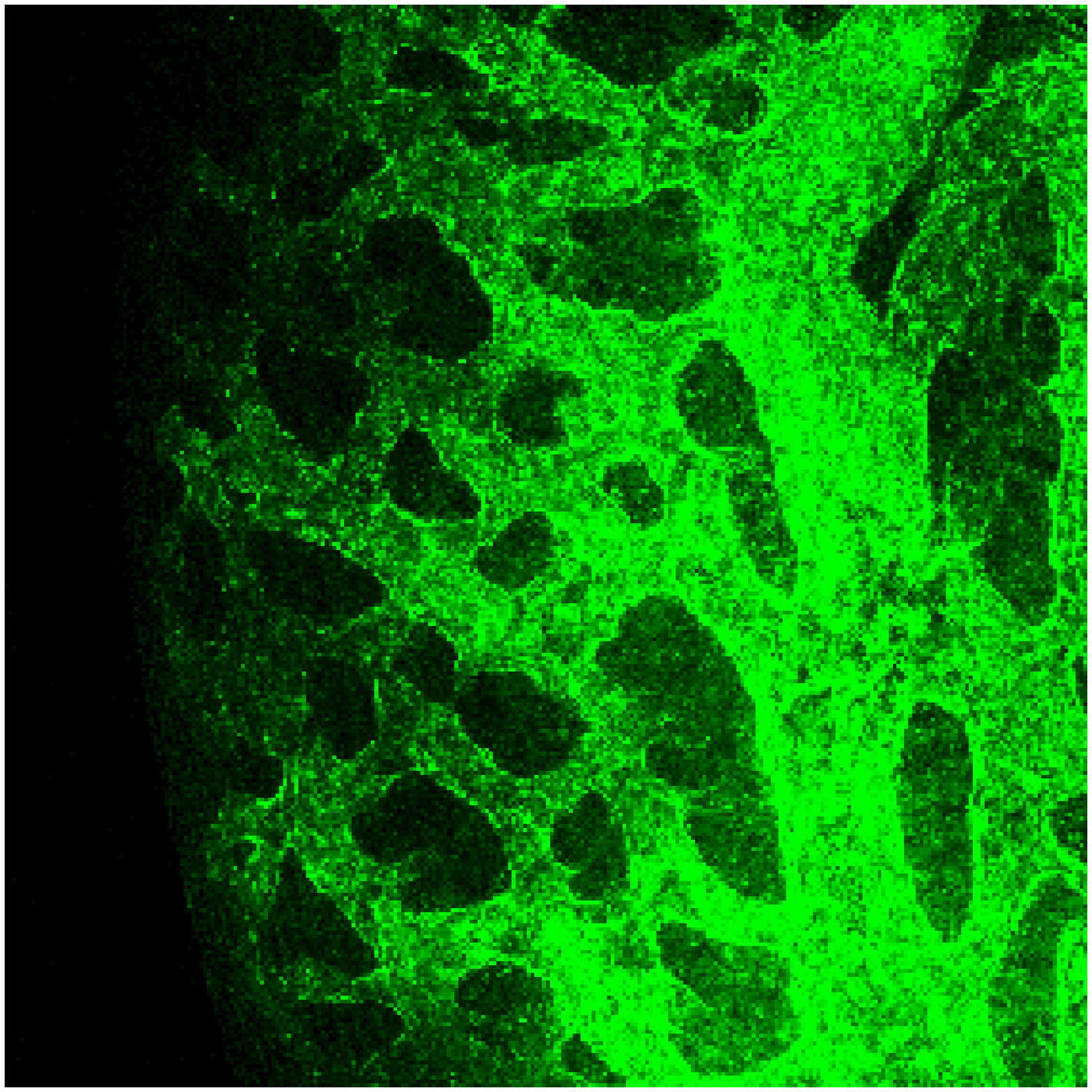}

\textsf{\textbf{C}} \includegraphics[%
  scale=0.3]{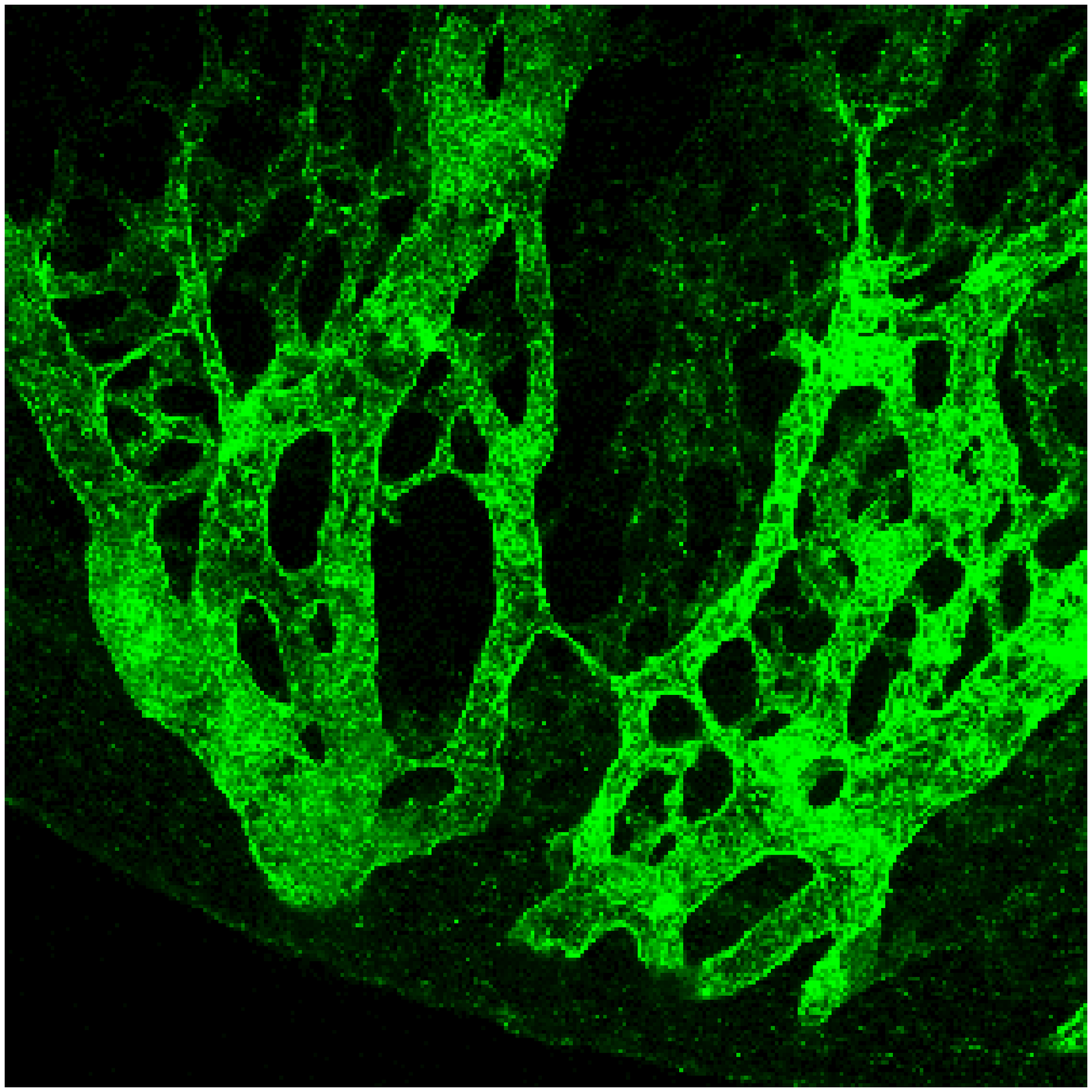} \textsf{\textbf{D}}\includegraphics[%
  scale=0.3]{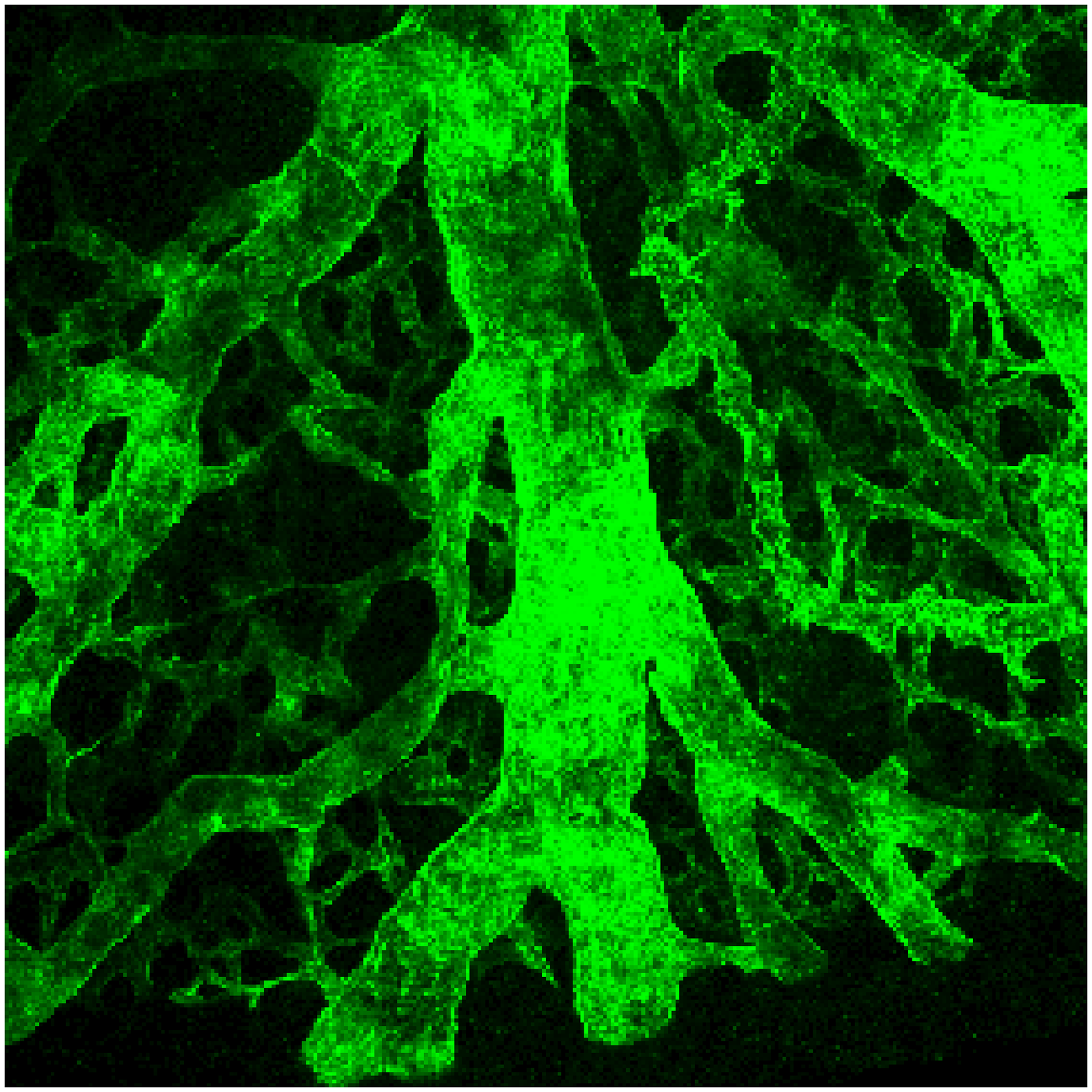}
\caption{Vascular networks formed by vasculogenesis in
chick embryo brain, at various stages of development, classified according
to Hamilton and Hamburger (HH). 
\textbf{A}: HH stage 17, corresponding to 52-64 hours; 
\textbf{B}: HH stage 20 (70-72 hours); 
\textbf{C,D}: HH stage 26 (5 days).}
\miolabel{fig:vascular}
\end{figure}

\section{Mathematical model of blood vessel growth}

The multidimensional Burgers' equation is a well-known paradigm in the
study of pattern formation. It gives a coarse grained hydrodynamic
description of the motion of independent agents performing rectilinear
motion and interacting only at very short ranges. These equations have
been utilized to describe the emergence of structured patterns in many
different physical settings (see e.g. \cite{SZ89,KPZ86}). In the early
stages of dynamics, each particle moves with a constant velocity,
given by a random statistical distribution. This motion gives rise to
intersection of trajectories and formation of shock waves. After the
birth of these local singularities regions of high density grow and
form a peculiar network-like structure. The main feature of this
structure is the existence of comparatively thin layers and filaments
of high density that separate large low-density regions.

In order to study and identify the factors influencing blood vessel
formation one has to take into account evidence suggesting that cells
do not behave as independent agents, but rather exchange information
in the form of soluble chemical factors. This leads to the model proposed
by Gamba et al. in ~\cite{GAC+03} and Serini et al. in~\cite{SAG+03}.
The model describes the motion of a fluid of randomly seeded independent
particles which communicate through emission and absorption of a soluble
factor and move toward its concentration gradients.%

\subsection{Model equations}

The cell population is described by a continuous density
$n(\mathbf{x},t)$, where $\mathbf{x}\in \mathbf{R}^{d}$ ($d=2,3$) is
the space variable, and $t\geq0$ is the time variable. The population
density moves with velocities $\mathbf{v}(\mathbf{x},t)$, that are
stimulated  by chemical gradients of a soluble factor. The
chemoattractant soluble factor is described by a scalar chemical
concentration field $c(\mathbf{x},t)$. It is supposed to be released 
by the cells, diffuse, and degrade in a finite time, in agreement with 
experimental observations.

The dynamics of the cell density can be described by coupling three
equations. 
The first one is the mass conservation law for cell matter, which expresses the conservation
of the number of cells. The second one is a momentum balance law that 
takes into account the phenomenological chemotactic force, the dissipation
by interaction with the substrate, the phenomenon of cell directional
persistency along their trajectories and a term implementing an
excluded volume constraint \cite{GAC+03,AGS04}. Finally there is a
reaction-diffusion equation for the production, degradation and
diffusion of the concentration of the  chemotactic factor.
One then has the following system: 
\begin{subequations}
\miolabel{Naldi:modello}
\begin{align}
&\frac{\partial n}{\partial t} + \nabla \cdot (n \mathbf{v}) = 0 \miolabel{Naldi:modello:n}\\
&\frac{\partial \mathbf{v}}{\partial t} + \mathbf{v} \cdot \nabla
\mathbf{v}
 = \mu(c) \nabla c - \nabla \phi(n) - \beta(c) \mathbf{v} \miolabel{Naldi:modello:v}\\
&\frac{\partial c}{\partial t} = D\mathrm{\Delta} c + \alpha(c) n - \frac{c}{\tau} \miolabel{Naldi:modello:c}
\end{align}
\end{subequations} 
where $\mu$ measures the cell response to the chemotactic factor, while
$D$ and $\tau$ are respectively the diffusion coefficient and the
characteristic degradation time of the soluble chemoattractant. The
function $\alpha$ determines the rate of release of the chemical
factor. The friction term $-\beta\mathbf{v}$ mimics the dissipative interaction of
the cells with the extracellular matrix. 

A simple model can be obtained
by assuming that the cell sensitivity $\mu$, the rate of release of
the chemoattractant $\alpha$ and the friction coefficient $\beta$ are
constant. A more realistic description may be obtained including
saturation effects as functional dependencies of the
aforementioned coefficients on the concentration $c$.

The term $\nabla \phi(n)$ is a density dependent pressure term, where
$\phi(n)$ is zero for low densities, and increases for densities above
a suitable threshold. This pressure is a phenomenological term which
models short range interaction between cells and the fact that cells
do not interpenetrate.

We observe that, at low density $n$ and for small chemoattractive
gradients, (\ref{Naldi:modello:v}b) is an inviscid Burgers' equation
for the velocity field $\mathbf{v}$ \cite{Bur74}, coupled to the
standard reaction-diffusion equation (\ref{Naldi:modello:c}c) and the
mass conservaton law (\ref{Naldi:modello:n}a). 

Since in the early stages of development almost all intraembryonic
mesodermal tissues contain migrating endothelial precursors, we use initial conditions
representing a randomly scattered distribution of cells, \emph{i.e.},
we throw an assigned number of cells in random positions inside the
cubic box, with zero initial velocities and zero initial concentration
of the soluble factor, with a single cell given initially by a
Gaussian bump of width $\sigma$ of the order of the average cell
radius ($\simeq15\mu\mathrm{m})$ and unitary weight in the integrated
cell density field $n$.

In order to model the fact that closely packed cells resist to
compression, a phenomenological, density dependent, pressure
$\nabla\phi(n)$ acting only when cells become close enough to each
other is introduced. The potential $\phi$ has to be monotonically
increasing and constant for $n<n_0$ where $n_0$ is the close-packing
density. Our simulations suggest that the exact functional form of
$\phi(n)$ is not relevant. For simplicity we choose
\begin{equation} \label{eq:pressione}
 \phi(n)=
   \begin{cases}
    B_{p}(n-n_{0})^{C_p} & n>n_{0}\\
    0 & n\leq n_{0}
   \end{cases}
\end{equation}

\subsection{Parameter values}

Fourier analysis of Eq. (\ref{Naldi:modello:c}) with constant
parameters and in the fast diffusion approximation $\partial
c/\partial t=0$ suggests that starting from the aformentioned
initial conditions, equation (\ref{Naldi:modello}) should develop
network patterns characterized by a typical length scale
$r_{0}=\sqrt{D\tau}$, which is the effective range of the interaction
mediated by soluble factors. As a matter of fact, Fourier components
$\hat c_{k}$ of the chemical field are related to the Fourier components of
the density field $\hat n_{k}$ by the relation
\[ \hat c_{k}=\frac{\alpha\tau \hat n_{k}}{D\tau k^{2}+1}.\]
This means that in equation (\ref{Naldi:modello}) wavelengths of the
field $n$ of order $r_{0}$ are amplified, while wavelengths
$\lambda\gg r_{0}$ or $\lambda\ll r_{0}$ are suppressed.

Initial conditions introduce in the problem a typical length scale
given by the average cell-cell distance $L/\sqrt{N}$, where $L$ is the
system size and $N$ the particle number. The dynamics, filtering
wavelengths \cite{TGL+06}, rearranges the matter and forms a network
characterized by the typical length scale $r_{0}$.

It is interesting to check the compatibility of the theoretical
prediction with physical data. From available experimental results
\cite{PNJ+99} it is known that the order of
magnitude of the diffusion coefficient for major angiogenic growth
factors is
$D=10^{-7}\hspace{0.25em}\mathrm{cm}^{2}\hspace{0.25em}\mathrm{s}^{-1}$.
In the experimental conditions that were considered in \cite{GAC+03} 
the half life of soluble factors is $64\pm7\hspace{0.25em}\mathrm{\min}$.
This gives $r_{0}\sim200\hspace{0.25em}\mu\mathrm{m}$, a value in good
agreement with experimental observations.

\subsection{Lower dimensional models}
In order to get some intuition about the typical system dynamics, 
we exploit the 1D version of model \eqref{Naldi:modello} to simulate
the ``collision'' of two cells. For small values of
$B_p$ and sufficiently high $C_p$ in (\ref{eq:pressione}), the two bumps merge into a single
one (see Fig. \ref{fig:press} left) which appears to be stationary,
as suggested also by the graphs
of the kinetic energy and of the momentum of inertia (Fig.
\ref{fig:press2} top). 
On the other hand a less smooth onset
of pressure obtained with larger $B_p$ or smaller $C_p$ leads to
forces overcoming the chemical attractive ones, making
the two bumps bounce back (Fig.
\ref{fig:press} right, Fig. \ref{fig:press2} bottom).  We observe that the better dynamics from the biological point of view is the first 
behavior with two bump coalescing.

\begin{figure}
\begin{center}
\begin{tabular}{c@{\hspace{15mm}}c}
\includegraphics[width=.45\textwidth]{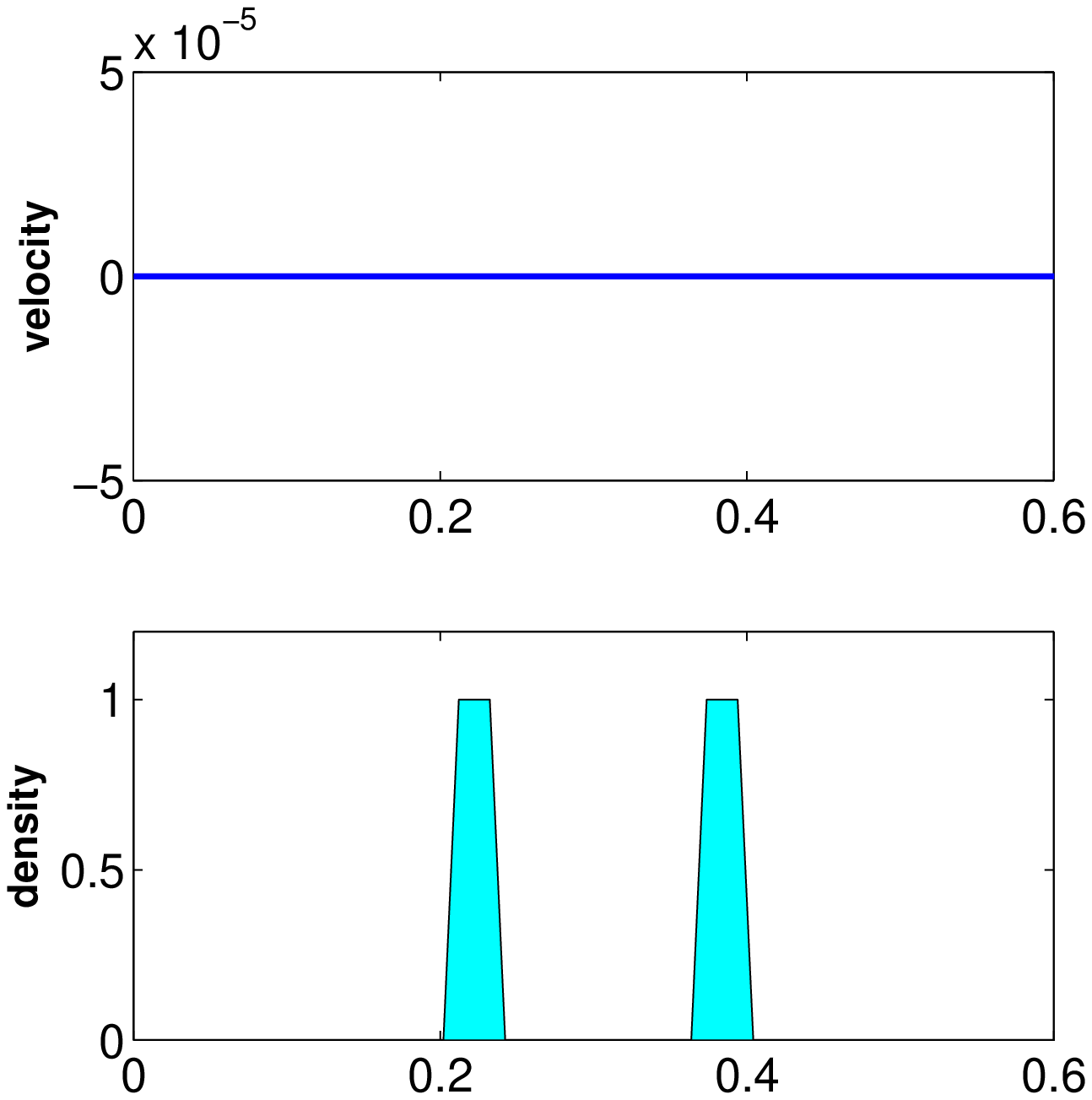}&
\includegraphics[width=.45\textwidth]{figure/rimbalzo_0.eps}\\[5mm]
\includegraphics[width=.45\textwidth]{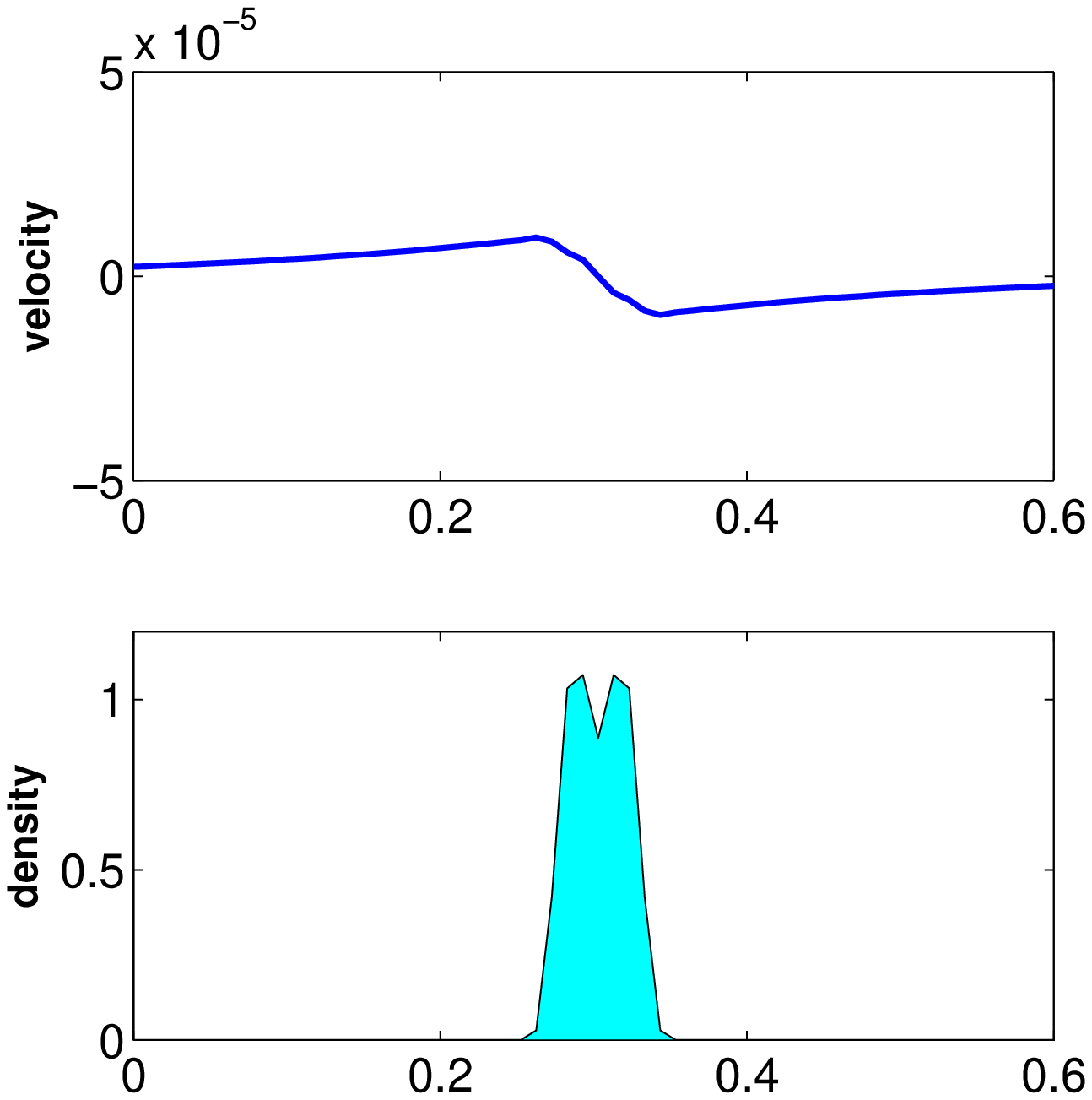}&
\includegraphics[width=.45\textwidth]{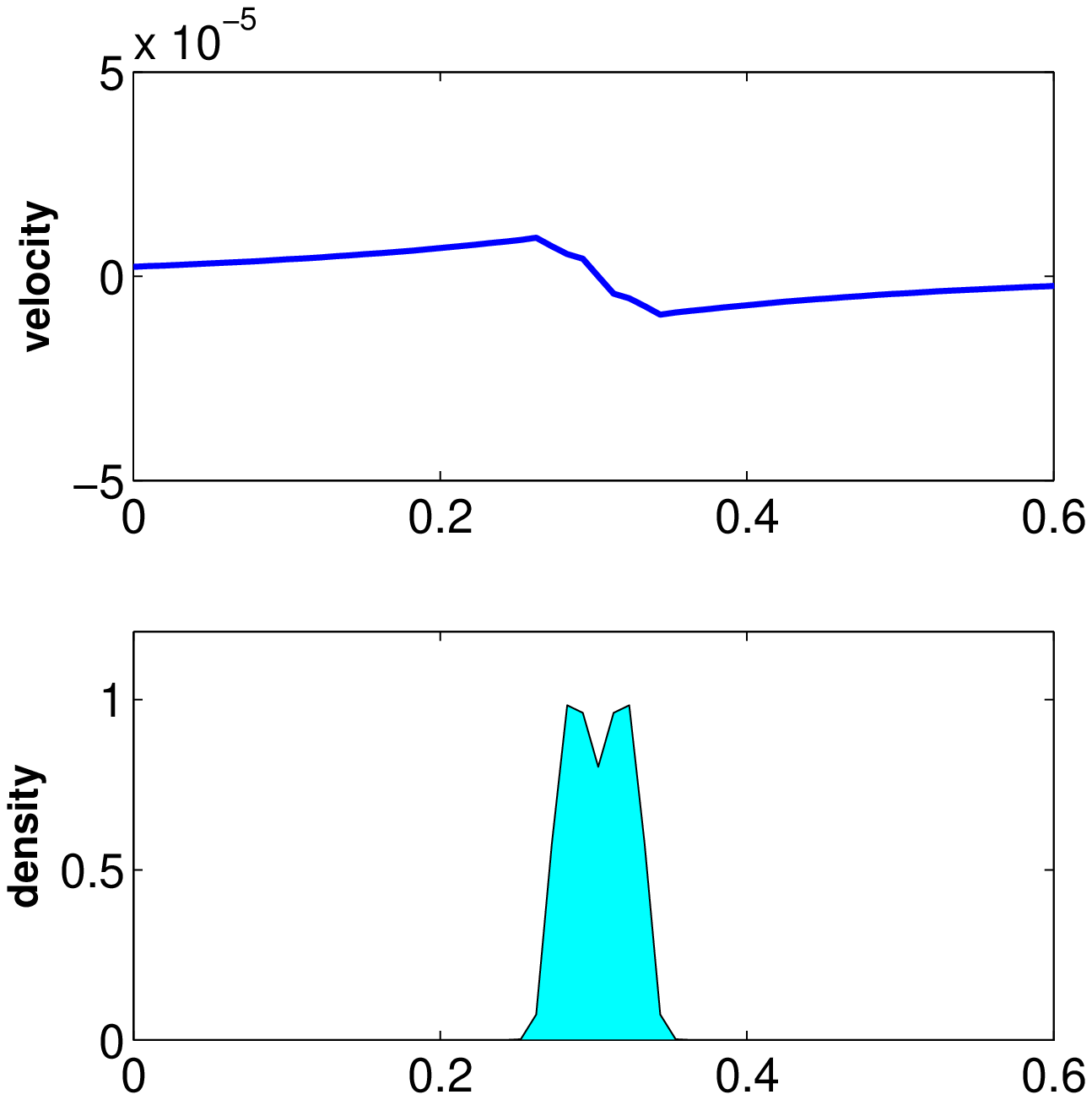}\\[5mm]
\includegraphics[width=.45\textwidth]{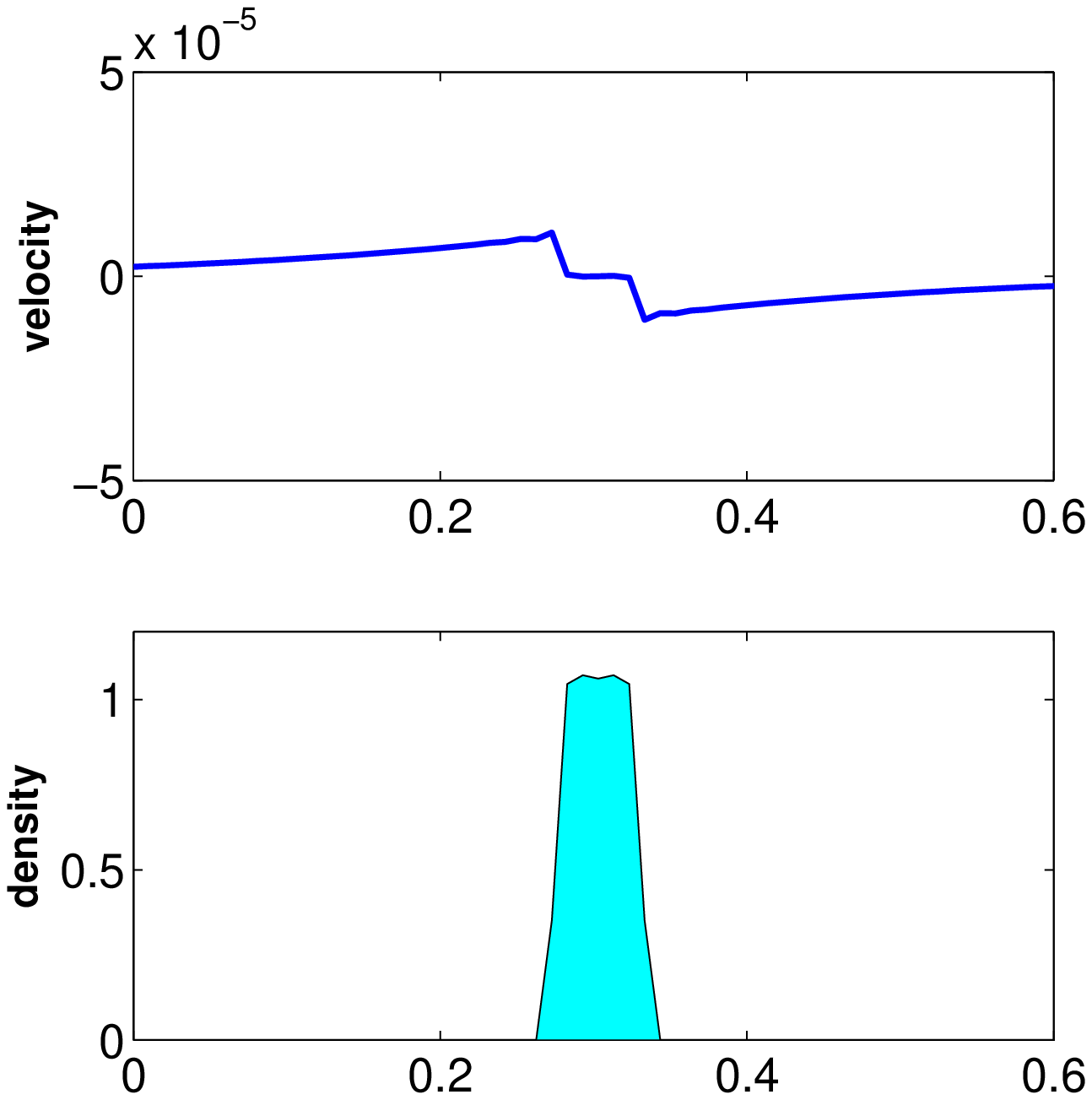}&
\includegraphics[width=.45\textwidth]{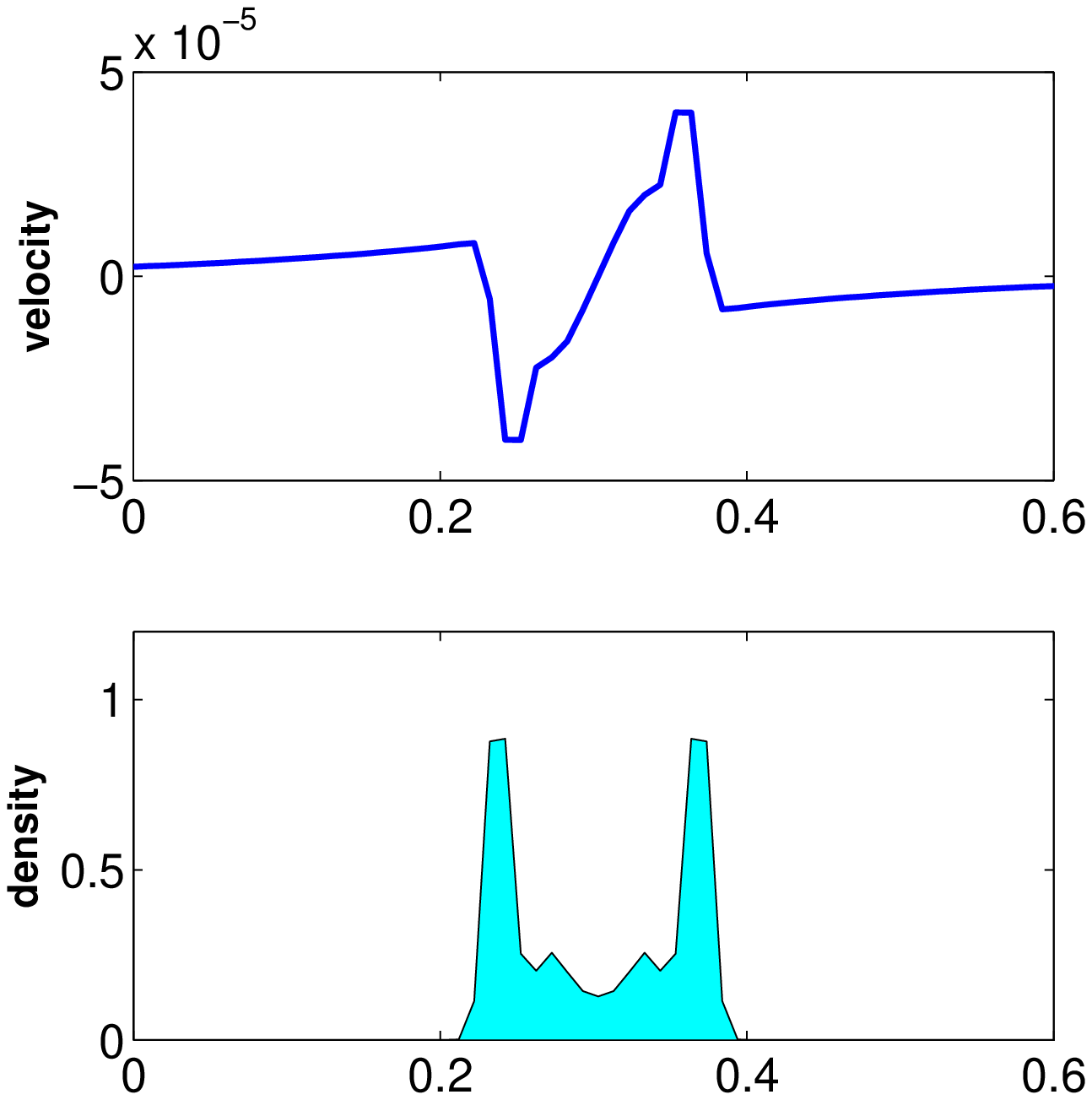}\\
\end{tabular}
\end{center}
\caption{Bump coalescence driven by chemotactic force and pressure.
In the first three rows the density and velocity fields at subsequent
instants of time are shown. In the last row we show the time evolution
of the kinetic energy and of the momentum of inertia.
Left column: $C_p=3$ and $B_p=10^{-3}$, leading to bump
coalescence.
Right column: $C_p=2$ and $B_p=10^{-1}$, leading to undesired
rebound of the two bumps.}
\miolabel{fig:press}
\end{figure}
\begin{figure}
\begin{center}
\begin{tabular}{c@{\hspace{15mm}}c}
\includegraphics[width=.9\textwidth]{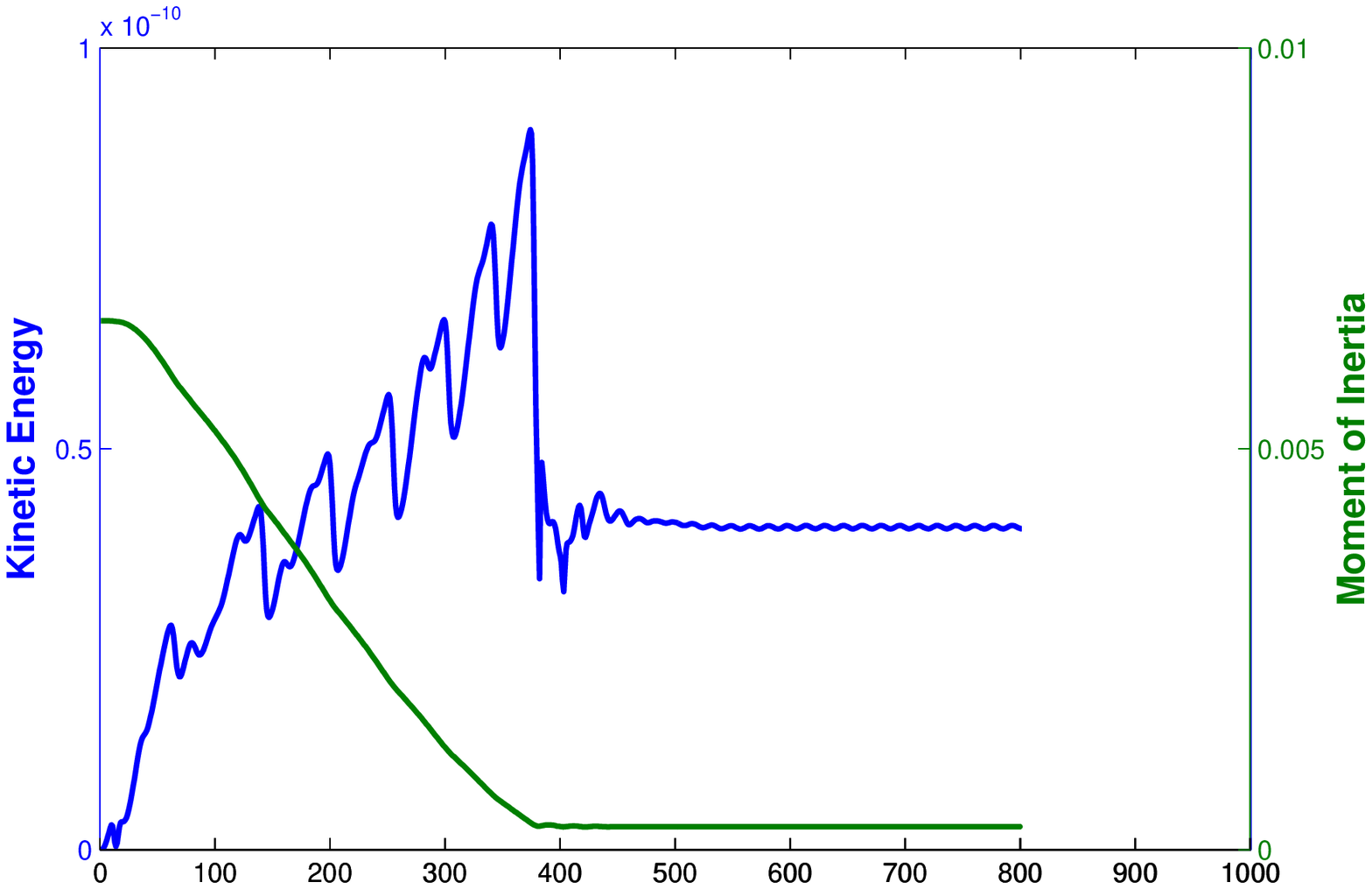}&\\[5mm]
\includegraphics[width=.9\textwidth]{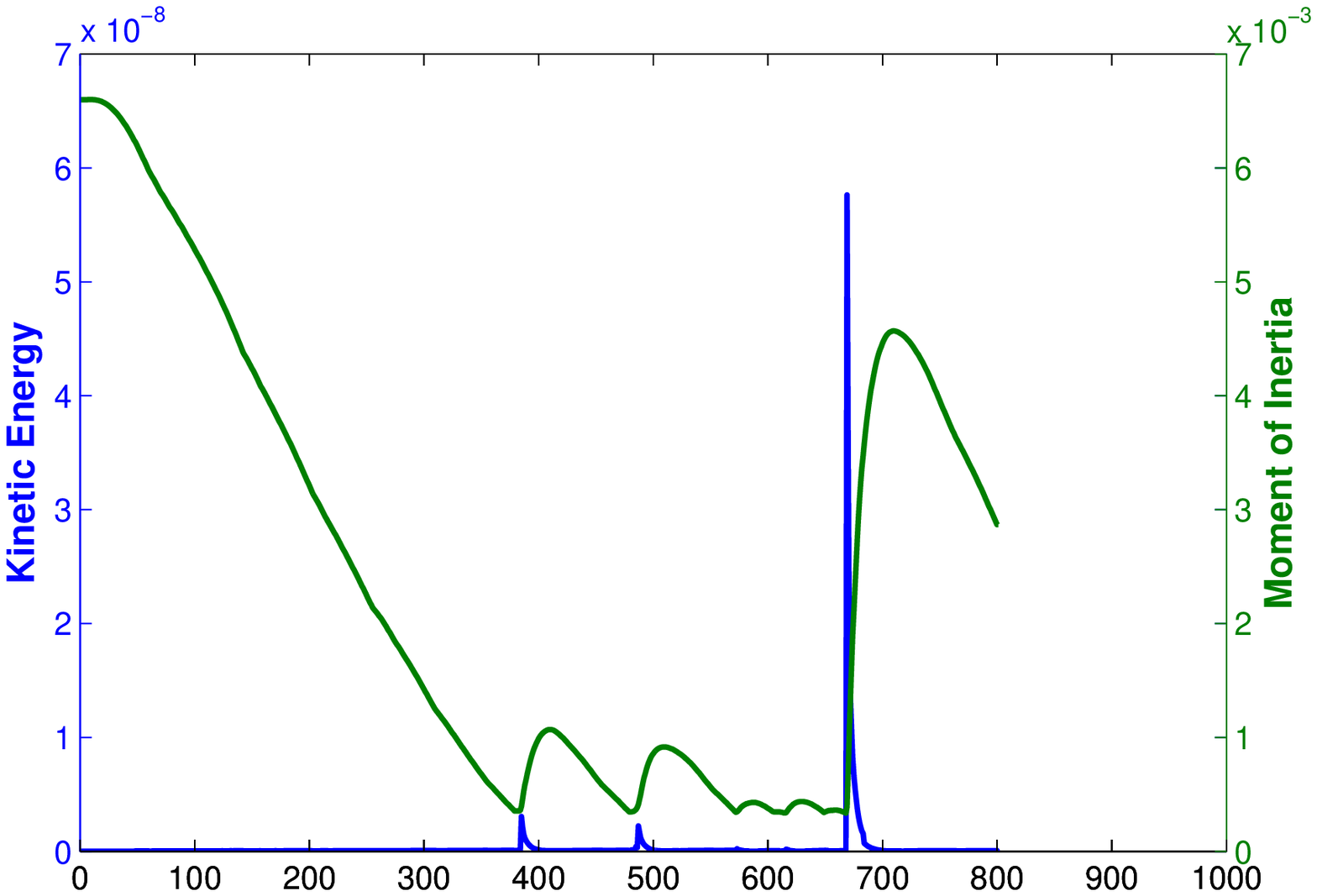}\\
\end{tabular}
\end{center}
\caption{Time evolution
of the kinetic energy and of the momentum of inertia.
Top: $C_p=3$ and $B_p=10^{-3}$, leading to bump
coalescence.
Bottom: $C_p=2$ and $B_p=10^{-1}$, leading to undesired
rebound of the two bumps.}
\miolabel{fig:press2}
\end{figure}

Biological observations suggest that the dynamics of cell changes
when they establish cell-cell contacts. It is reasonable to suppose
that a different genetic program is activated at this moment, disabling
cell motility. We therefore switch off cell motility as soon as the
cell concentration, signalled by chemoattractant emission, reaches
a given threshold. In this way the computational system is guaranteed
to reach a stationary state.

These effects can be taken into account using a non-constant
sensitivity $\mu(c)$, a non-linear emission rate $\alpha(c)$, or a
variable friction coefficient $\beta(c)$. We choose a threshold
$c_{0}$ and functions of the form
\begin{subequations} \miolabel{Naldi:mube}
\begin{align}
 \mu(c) &= \mu_0 [1- \tanh(c-c_0)] \miolabel{Naldi:mu} \\   
 \alpha(c) &= \alpha_0 [1- \tanh(c-c_0)]     \miolabel{Naldi:be}\\
 \beta(c) &= \beta_0 [1+ \tanh(c-c_0)] \miolabel{Naldi:ro}
\end{align}
\end{subequations}
The effect of the first two terms is that the sensitivity of the cells and 
their chemoattractant production is strongly damped when the concentration
$c$ reaches the threshold $c_{0}$. We did not observe a significant
dependence on the exact form of the damping function, provided that
it approximates a step function that is nonzero only when $c<c_{0}$.

$\beta(c)$, on the other hand has the effect of turning on a strong
friction term at locations of high chemoattractant concentration. 
We performed several tests and observed that the different
choices are approximately equivalent in freezing the system into a
network-like stationary state.

\section{Numerical methods}

Our scheme is based on a suitable relaxation approximation
\cite{JX95} of the mass conservation law  (\ref{Naldi:modello:n}) and the multidimensional Burgers
equation (\ref{Naldi:modello:v}) coupled with a  second order
finite-differences method for the reaction-diffusion equation (\ref{Naldi:modello:c}) of the chemotactic factor.
We point out that also for the  last equation (\ref{Naldi:modello:c})  we
could consider a relaxation approximation \cite{JPT99,NP00}
in order to deal with the system (\ref{Naldi:modello}) in an uniform way,
but we prefer to adopt here a simpler approach.

We first briefly review an extension of the approach proposed by Jin and
Xin in \cite{JX95} for a scalar conservation law to the case when a
source term is present
\begin{equation} \label{Naldi:conslaw}
\frac{\partial u}{\partial t}+\frac{\partial}{\partial x}f(u)=g(u). 
\end{equation}
Introducing an auxiliary variable $j$ that plays the role of a
physical flux we consider the following relaxation system:
\begin{subequations}
\label{Naldi:relaxation}
\begin{align}
&\frac{\partial u}{\partial t} + \frac{\partial j}{\partial x} = g(u)\\
&\frac{\partial j}{\partial t} + a \frac{\partial u}{\partial x} =
  -\frac1{\epsilon}(j-f(u)) ,
\end{align}
\end{subequations}
where $\epsilon$ is a small positive parameter, called relaxation
time, and $a$ is a suitable positive constant.  Formally,
Chapman-Enskog expansion justifies the agreement of the solutions of
the relaxation system with the solutions of the equation
\begin{equation}
\frac{\partial u}{\partial t}+\frac{\partial}{\partial
  x}f(u)=g(u)+\epsilon\frac{\partial}{\partial
  x}\left((a-f'(u)^{2})\frac{\partial u}{\partial
  x}\right),\miolabel{Naldi:Chapman}
\end{equation}
which is a first order approximation of the original balance law
\eqref{Naldi:conslaw}.

It is also clear that \eqref{Naldi:Chapman} is dissipative, provided
that the subcharacteristic condition $a>f'(u)^{2}$ is satisfied. We
would expect that appropriate numerical discretization of the
relaxation system \eqref{Naldi:relaxation} yields accurate
approximation to the original equation \eqref{Naldi:conslaw} when the
relaxation parameter $\epsilon$ is sufficiently small.

In view of its numerical approximation, the main advantage of the
relaxation system \eqref{Naldi:relaxation} over the original equation
\eqref{Naldi:conslaw} lies in the linear structure of the characteristic
fields and in the localized low order term and this avoids the use of time 
consuming Riemann solvers. Moreover, proper implicit time discretization 
can be exploited to overcome the stability constraints due to the stiffness and 
to avoid the use of non-linear solvers.

\newcommand{\pder}[2]{\frac{\partial{#1}}{\partial{#2}}}

We observe that system \eqref{Naldi:relaxation} is in the form
\begin{equation} \label{eq:z}
 \pder{z}{t} + \mathrm{div}f(z) = g(z) + \frac{1}{\epsilon}h(z)
\end{equation}
where $z=(u,j)^T$, $f(z)=(j,au)^T$, $g(z)=(g(u),0)^T$ and 
$h(z)=(0,j-f(u))^T$. When $\varepsilon$ is small, the presence of both
non-stiff and stiff terms, suggests the use of IMEX schemes
\cite{ARS97,KC03,PR05}.

Assume for simplicity to adopt a uniform time step $\Delta t$ and
denote with $z^n$ the numerical approximation at time $t_n=n\Delta t$,
for $n=0,1,\ldots$ In our case a $\nu$-stages IMEX scheme reads
\[ z^{n+1} = z^n 
     - \Delta t \sum_{i=1}^{\nu}
           \tilde{b}_i  \left[\pder{f}{x}(z^{(i)}) +g(z^{(i)})\right]
     + \frac{\Delta t}{\varepsilon} \sum_{i=1}^{\nu} b_i h(z^{(i)})
\]
where the stage values are computed as
\[z^{(i)} = z^n 
        -\Delta t \sum_{k=1}^{i-1}\tilde{a}_{i,k} \left[\pder{f}{x}(z^{(k)}) +g(z^{(k)})\right]
        + \frac{\Delta t}{\varepsilon} \sum_{k=1}^{i} a_{i,k} h(z^{(k)})
\]
Here $({a}_{ik},{b}_i)$ and  $(\tilde{a}_{ik},\tilde{b}_i)$ are a pair
of Butcher's tableaux of, respectively, a diagonally implicit and an
explicit Runge-Kutta schemes.

In this work we use the so-called relaxed schemes, that are obtained
letting $\varepsilon\rightarrow0$ in the numerical scheme for
\eqref{eq:z}. For these the first stage 
\[ 
 \left[\begin{array}{l} 
    {u^{(1)}}\\
    {j^{(1)}}
 \end{array}\right]
 = 
 \left[\begin{array}{l} u^{n}\\j^{n}\end{array}\right]
 +
 \frac{\Delta t}{\varepsilon} a_{1,1} 
   h\left(\left[\begin{array}{l}u^{(1)}\\j^{(1)}\end{array}\right]\right)
\]
becomes
\[ {u^{(1)}=u^n} \qquad {j^{(1)}=f(u^{(1)})},\]
then it reduces to $h(z^{(1)})=0$.
While the second stage, $i=2$, reads 
\[z^{(2)} = z^n 
    - \Delta{t}\tilde{a}_{2,1}{\left[\pder{f}{x}(z^{(1)}) + g(z^{(1)})\right]} 
    + \frac{\Delta t}{\varepsilon} a_{2,1}
           \underbrace{{h(z^{(1)})}}_{{\equiv 0}}
    + \frac{\Delta t}{\varepsilon} a_{2,2} {h(z^{(2)})}
\]
which implies that $h(z^{(2)})=0$.

Summarizing, the relaxed scheme yields an alternation of
relaxation steps
\[ h(z^{(i)})=0 \qquad 
   \text{ i.e. }\, 
     {j^{(i)}=f(u^{(i)})}
\]
and transport steps where we advance for time $\tilde{a}_{i,k}\Delta{t}$
\[ \pder{z}{t} + \mathrm{div} f(z) = g(z) \]
with initial data $z=z^{(i)}$ retain only the first component and
assign it to $u^{(i+1)}$.

Finally the value of $u^{n+1}$ is computed as \(u^n+\sum \tilde{b}_i u^{(i)}\).

In order to obtain a relaxation approximation of the first and second
equation of \eqref{Naldi:modello} we
rewrite them in
conservative form, introducing the moment
$\mathbf{p}(\mathbf{x},t)=n(\mathbf{x},t)\mathbf{v}(\mathbf{x},t)$:
\begin{subequations}
\label{modello:cons}
\begin{align}
&\frac{\partial n}{\partial t} + \nabla \cdot \mathbf{p} = 0 \miolabel{modello:cons:n}\\
&\frac{\partial \mathbf{p}}{\partial t} + \nabla\cdot\left(\mathbf{n}\mathbf{v}\otimes\mathbf{v}\right)
 = n\mu \nabla c - n\nabla \phi(n) - \beta \mathbf{p} \miolabel{modello:cons:p}
\end{align}
\end{subequations}

Introducing the variable $\mathbf{u}=(n,\mathbf{p})^T$ and the
auxiliary flux $\mathbf{w}$, the relaxation system reads
\begin{subequations}
\label{relax:cons}
\begin{align}&\frac{\partial \mathbf{u}}{\partial t} + \nabla \cdot \mathbf{w} = G(\mathbf{u},\mathbf{w},c) \miolabel{relax:cons:u}\\
&\frac{\partial \mathbf{w}}{\partial t} + A\nabla \cdot \mathbf{u}= -\frac{1}{\varepsilon}\left(\mathbf{w}-F(\mathbf{u})\right)
\end{align}
\end{subequations}
where $G(\mathbf{u},\mathbf{w},c)=(0,n\mu \nabla c - n\nabla \phi(n) -
\beta \mathbf{p})^T$,
$F(\mathbf{u})=(\mathbf{p},\mathbf{n}\mathbf{v}\otimes\mathbf{v})$ and $A$ is a
suitable diagonal matrix whose positive diagonal elements verify a
subcharacteristic condition.  
As we previously remarked, our relaxed scheme takes alternatively an implicit step and an explicit one: the explicit step involves the
computation of the flux $\nabla \cdot \mathbf{w}$ and the evaluation
of the non stiff source term $G$. In particular we compute $\nabla c$
and $\nabla \phi(n)$ using a second order difference scheme. 
 
In the following we describe for simplicity the fully discrete scheme in one dimensional case. We
introduce the spatial grid points $x_{j}$ with uniform mesh width
$h=x_{j+1}-x_{j}$. As usual, we denote by $u_j^n$ the
approximate cell average of a quantity $u$ in the cell
$[x_{j-1/2},x_{j+1/2}]$ at time $t_n$ and by $u_{j+1/2}^n$ the
approximate point value of $u$ at $x=x_{j+1/2}$ and $t=t_n$.
A spatial discretization to (\ref{relax:cons}) in conservation form
can be written as 
\begin{subequations}
\label{relax:consnum}
\begin{align}&\frac{\partial \mathbf{u}_j}{\partial t} + \frac{1}{h}\left(\mathbf{w}_{j+1/2}-\mathbf{w}_{j-1/2}\right) = G(\mathbf{u}_j,\mathbf{w}_j,c_j) \miolabel{relax:cons:unum}\\
&\frac{\partial \mathbf{w}_j}{\partial t} + \frac{1}{h}A\left(\mathbf{u}_{j+1/2}-\mathbf{u}_{j-1/2}\right)= -\frac{1}{\varepsilon}\left(\mathbf{w}_j-F(\mathbf{u}_j)\right).
\end{align}
\end{subequations}
In order to compute the numerical fluxes $\mathbf{w}_{j\pm1/2}$, 
we consider the characteristic variables $\mathbf{w}\pm
A^{1/2}\mathbf{u}$ that travel with constant velocities $\pm A^{1/2}$, and so the semidiscrete system becomes 
diagonal. Now we have to apply a numerical approximation to $\mathbf{w}\pm
A^{1/2}\mathbf{u}$. A first idea is to apply a ENO or WENO approach (see e.g. \cite{Shu97}), to build an high order reconstruction, coupled with a suitable IMEX scheme. The drawback is the high computational costs, especially in a multidimensional framework.
Therefore we chose a suitable compromise between the
computational cost and the accuracy, using a second order TVD scheme. The
numerical flux that we use is obtained coupling an upwind scheme and the
Lax-Wendroff method by a non linear flux limiter \cite{leveque90}. 
Namely the high order flux $F(U)$ for a generic variable $U$ consists of the low order term $F_L(U)$
plus a second order correction $F_H(U)$:
\[
F(U)=F_L(U)+\Psi(U)(F_H(U)-F_L(U)) 
\]
where $\Psi$ is the flux limiter. When the data $U$ is smooth, then
$\Psi(U)$ should be near $1$, while near a discontinuity we want
$\Psi(U)$ close to $0$. The idea consists in the selection of a high
order flux $F_H$ that works well in smooth regions and of a low order
flux $F_L$ which behaves well near discontinuities.   

In our schemes we considered the upwind scheme as a low order
flux for the characteristic variables
\[
F_L((\mathbf{w}+A^{1/2}\mathbf{u})_{j+1/2})=(\mathbf{w}+A^{1/2}\mathbf{u})_j, \quad 
F_L((\mathbf{w}-A^{1/2}\mathbf{u})_{j+1/2})=(\mathbf{w}-A^{1/2}\mathbf{u})_{j+1}
\]
and the Lax-Wendroff scheme as a high order flux for the same
variables 
\[
\begin{array}{ll}
F_H((\mathbf{w}\pm
A^{1/2}\mathbf{u})_{j+1/2})=&\frac{A^{1/2}}{2}((\mathbf{w}\pm
A^{1/2}\mathbf{u})_{j+1}+(\mathbf{w}\pm 
A^{1/2}\mathbf{u})_{j})\\
&-\frac{\lambda A^{1/2}}{2}((\mathbf{w}\pm
A^{1/2}\mathbf{u})_{j+1}-(\mathbf{w}\pm A^{1/2}\mathbf{u})_{j})
\end{array}
\]
where $\lambda=\Delta t/h$ (we advance of one time step).

Letting 
\[
\mathbf{\Theta}_j^{\pm}=\left( \frac{(\mathbf{w}\pm
A^{1/2}\mathbf{u})_{j}^n-(\mathbf{w}\pm
A^{1/2}\mathbf{u})_{j-1}^n}{(\mathbf{w}\pm
A^{1/2}\mathbf{u})_{j+1}^n-(\mathbf{w}\pm
A^{1/2}\mathbf{u})_{j}^n} \right)^{\pm1}\,,
\]
the fully discrete scheme for the variable $\mathbf{u}$ using Euler
method to advance in time is the following 
\[
\begin{array}{ll}
\mathbf{u}_j^{n+1}=&\mathbf{u}_j^{n}+\frac{\lambda A^{1/2}}{2}
(\mathbf{u}_{j+1}^{n}-2\mathbf{u}_j^{n}+\mathbf{u}_{j-1}^{n})-\frac{\lambda}{2}(\mathbf{w}_{j+1}^{n}-\mathbf{w}_{j-1}^{n})\\
&\Delta t\frac{\mathbb{I}-\lambda A^{1/2}}{4}(-s_j^++s_{j-1}^+ +s_{j+1}^- - s_{j}^-),
\end{array}
\]
with
\begin{equation}
\label{eq:s}
s_j^{\pm}=\frac{1}{h}(\pm
A^{1/2}\mathbf{u}_{j\pm1}^n+\mathbf{w}_{j\pm1}^n\mp A^{1/2}\mathbf{u}_{j}^n-\mathbf{w}_{j}^n)\Psi(\Theta_j^{\pm}).
\end{equation}
After the substitution of the relaxing step we get
\[
\label{eq:relaxed}
\begin{array}{ll}
\displaystyle
\mathbf{u}_j^{n+1}=&\mathbf{u}_j^{n}+\frac{\lambda A^{1/2}}{2}
(\mathbf{u}_{j+1}^{n}-2\mathbf{u}_j^{n}+\mathbf{u}_{j-1}^{n})-\frac{\lambda}{2}(F(\mathbf{u}_{j+1}^{n})-F(\mathbf{u}_{j+1}^{n}))\\
\displaystyle
&\Delta t\frac{\mathbb{I}-\lambda A^{1/2}}{4}(-s_j^+ +s_{j-1}^+ +s_{j+1}^- - s_{j}^-),
\end{array}
\]
where $s^{\pm}$ is obtained from (\ref{eq:s}) letting
$\mathbf{w}=F(\mathbf{u})$.
The scheme can be put in a conservative form and it is possible to
prove its consistency by standard technique \cite{leveque90}.
In order to prove a TVD stability, we write 
\begin{equation}
\label{eq:harten}
\begin{array}{ll}
\displaystyle
\mathbf{u}_{j+1}^{n+1}-\mathbf{u}_{j}^{n+1}=&(1-\mathbf{C}_{j}^n-\mathbf{D}_{j}^n)(\mathbf{u}_{j+1}^{n}-\mathbf{u}_{j}^{n})+\mathbf{C}_{j-1}^n(\mathbf{u}_{j}^{n}-\mathbf{u}_{j-1}^{n})\\
\displaystyle
&+\mathbf{D}_{j+1}^n(\mathbf{u}_{j+2}^{n}-\mathbf{u}_{j+1}^{n})+E_{j+1/2},
\end{array}
\end{equation}
where
\[
\begin{array}{lll}
\displaystyle
\mathbf{C}_j^n&=&\frac{\lambda}{2}\left(A^{1/2}+\frac{F(\mathbf{u}_{j+1}^{n})-F(\mathbf{u}_{j}^{n})}{\mathbf{u}_{j+1}^{n}-\mathbf{u}_{j}^{n}}\right)\\
\displaystyle
\mathbf{D}_j^n&=&\frac{\lambda}{2}\left(A^{1/2}-\frac{F(\mathbf{u}_{j+1}^{n})-F(\mathbf{u}_{j}^{n})}{\mathbf{u}_{j+1}^{n}-\mathbf{u}_{j}^{n}}\right)\\
\displaystyle
E_{j+1/2}^n&=&\Delta t\frac{1-\lambda A^{1/2}}{4}\left(s_{j+2}^--2s_{j+1}^-+s_{j}^--s_{j+1}^++2s_{j}^+-s_{j-1}^+\right)
\end{array}
\]
where we notice that $\mathbf{C}$ and $\mathbf{D}$ are non negative.

The coefficient $E$ can be written in terms of $\mathbf{C}$ and $\mathbf{D}$, in fact
\[
s_j^+=\frac{2}{\lambda\Delta x}\mathbf{C}_i^n\Psi(\mathbf{\Theta}_i^+)(\mathbf{u}_{j+1}^{n}-\mathbf{u}_{j}^{n})), \quad s_j^-=-\frac{2}{\lambda\Delta x}\mathbf{D}_{i-1}^n\Psi(\mathbf{\Theta}_i^-)(\mathbf{u}_{j}^{n}-\mathbf{u}_{j-1}^{n})).
\]
We can rewrite (\ref{eq:harten}) in the following form
\begin{equation}
\label{eq:tvd}
\begin{array}{ll}
\displaystyle
\mathbf{u}_{j+1}^{n+1}-\mathbf{u}_{j}^{n+1}=&(\mathbf{u}_{j+1}^{n}-\mathbf{u}_{j}^{n})\left[(1-\mathbf{C}_{j}^n-\mathbf{D}_{j}^n)+(1-\lambda A^{1/2})(\mathbf{D}_j^n\Psi_{j+1}^-+\mathbf{C}_j^n\Psi_{j}^+)\right] \\
&+(\mathbf{u}_{j}^{n}-\mathbf{u}_{j-1}^{n})\left[\mathbf{C}_{j-1}^n-\frac{1-\lambda A^{1/2}}{2}(\mathbf{D}_{j-1}^n\Psi_{j}^-+\mathbf{C}_{j-1}^n\Psi_{j-1}^+)\right]\\
&+(\mathbf{u}_{j+2}^{n}-\mathbf{u}_{j+1}^{n})\left[\mathbf{D}_{j+1}^n-\frac{1-\lambda A^{1/2}}{2}(\mathbf{D}_{j+1}^n\Psi_{j+2}^-+\mathbf{C}_{j+1}^n\Psi_{j+1}^+)\right]
\end{array}
\end{equation}
It's easy to see that under the CFL condition $\|\lambda\sqrt{\max\{a_i\}}\|\leq1$, where $a_i$ are the positive diagonal
elements of the matrix $A$, and using the fact that the flux limiter verifies
\[
0\leq\frac{\Psi(\mathbf{\Theta})}{\mathbf{\Theta}}\leq2, \qquad 0\leq\Psi(\mathbf{\Theta})\leq2,
\]
we have
\[
\begin{array}{rrr}
(1-\mathbf{C}_{j}^n-\mathbf{D}_{j}^n)+(1-\lambda A^{1/2})(\mathbf{D}_j^n\Psi_{j+1}^-+\mathbf{C}_j^n\Psi_{j}^+)&\geq&0\\
\mathbf{C}_{j-1}^n-\frac{1-\lambda A^{1/2}}{2}(\mathbf{D}_{j-1}^n\Psi_{j}^-+\mathbf{C}_{j-1}^n\Psi_{j-1}^+)&\geq&0\\
\mathbf{D}_{j+1}^n-\frac{1-\lambda A^{1/2}}{2}(\mathbf{D}_{j+1}^n\Psi_{j+2}^-+\mathbf{C}_{j+1}^n\Psi_{j+1}^+)&\geq&0
\end{array}
\]
and so we can deduce that our scheme is TVD stable from Harten's Theorem \cite{Har84}.

In the case of multidimensions, a similar discretization can be
applied to each space dimension \cite{JX95,JPT99,NP00}.  Then, since
the structure of the multidimensional relaxation system is
similar to the 1D system, the numerical implementation for higher
dimensional problems, based on additive dimensional splitting, is not
much harder than for 1D problems.

 \begin{figure}
 \begin{center}
 \includegraphics[scale=0.5]{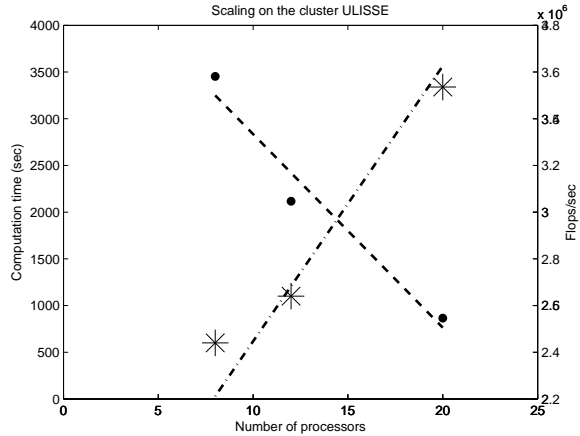}
 \end{center}
 \caption{Scaling of the 3D algorithm on the ULISSE cluster. Dots
   represent execution time (s) and asterisks the number of
   Mflops/s for our numerical algorithm. Dashed and dash-dot lines
   are linear interpolations}
 \miolabel{fig:scaling}
 \end{figure}
For our threedimensional problem the computational cost is quite high and can be reduced using
parallel computing: the semilinearity of relaxation systems, together
with our suitably chosen discretizations, provides parallel algorithm
with almost optimal scaling properties. 
In particular the domain is divided in smaller
subdomains and each subdomain is assigned to a processor. The computations of all non linear terms involve only
pointwise evaluations and it is easy to perform these tasks in a local
way. Only point near the interfaces between different subdomain need
to be communicated in the transport step. 
We implemented these algorithm
on a high performance cluster for parallel computation installed at
the Department of Mathematics of the University of Milano
(http://cluster.mat.unimi.it/). 
The scaling properties of the algorithm
are shown in Fig. \ref{fig:scaling} and are essentially due to the
exclusive use of matrix-vectors operations and to the avoidance of
solvers for linear or non-linear systems.

\section{Numerical results}

We perform threedimensional numerical simulations of model
\eqref{Naldi:modello} on a cubic box with side of length
$L=1\mathrm{mm}$, with periodic boundary conditions. The initial
condition is assigned in the form of a set of gaussian bumps with
$\sigma=15\mathrm{\mu{m}}$ scattered in the cube with uniform
probability and having zero initial velocity.

Biochemical data \cite{SAG+03} suggest the values 
$D=10^{-3}\mathrm{mm^{2}/s}$ and $\tau=4000~\mathrm{s}$ for the
diffusion constant and the chemoattractant decay rate. We fix
the other constant parameters by dimensional analysis and fitting to the
characteristic scales of the biological system. In
particular, we choose: $\mu_{0}=10^{-11}\mathrm{mm}^{4}/\mathrm{s}^{3}$,
$\alpha=1\mathrm{s}^{-1}$, $\beta=10^{-3}\mathrm{s}^{-1}$.
For the coefficients in the expression  \eqref{eq:pressione} of the pressure function $\phi$
we take $n_0=1.0$,$C_p=3$ and $B_p=10^{-3}$.

Very fine grids have to be used in order to resolve the details of
the $n(\mathbf{x},t)$ field, which may contain hundreds of small bumps,
each representing a single cell. Since
each cell has radius $\sigma=15\mathrm{\mu m}$, one needs a grid
spacing such that $\mathrm{\Delta{x}}<10\mathrm{\mu m}$ and therefore
grids of at least $100^3$ cells for a cubic domain of $1\mathrm{mm}$ side.

\begin{figure}
\includegraphics[scale=0.6,angle=-90]{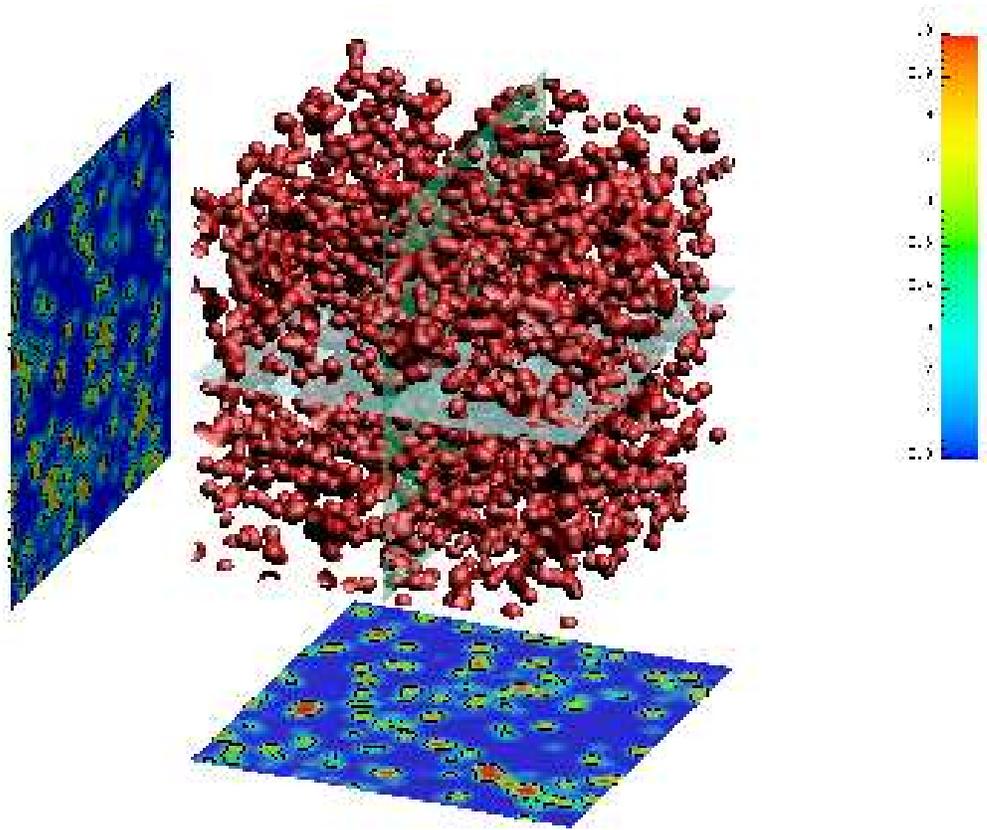}
\caption{Initial state of a numerical simulation with
  $2500\mathrm{ cells/mm}^3$. The colorbar on the right is referred to
  the coloring of the cross sections. The red three-dimensional
  isosurface corresponds to the black contour lines in the cross
  sections}
\miolabel{fig:3d:in}
\end{figure}

\begin{figure}
\includegraphics[scale=0.6,angle=-90]{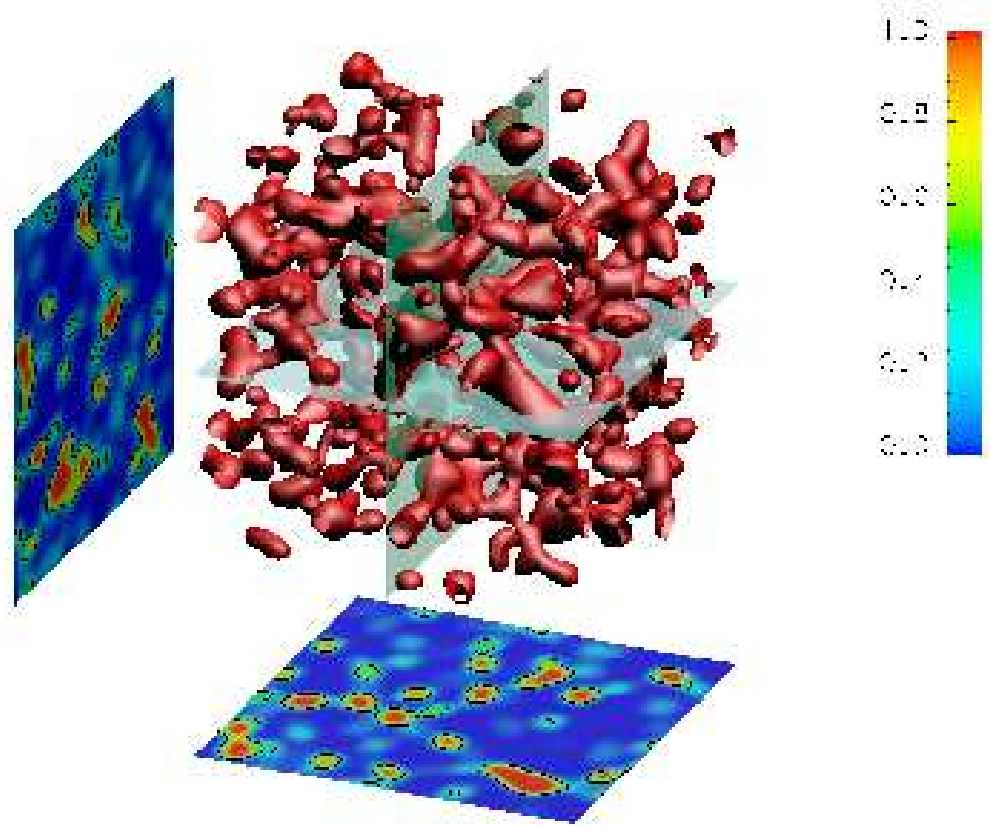}
\caption{Transient state of the evolution of the initial state
  depicted in Fig. \ref{fig:3d:in} according to model \eqref{Naldi:modello}. The
  initial formation of network-like structures is observed}
\miolabel{fig:3d:med}
\end{figure}

\begin{figure}
\includegraphics[scale=0.6,angle=-90]{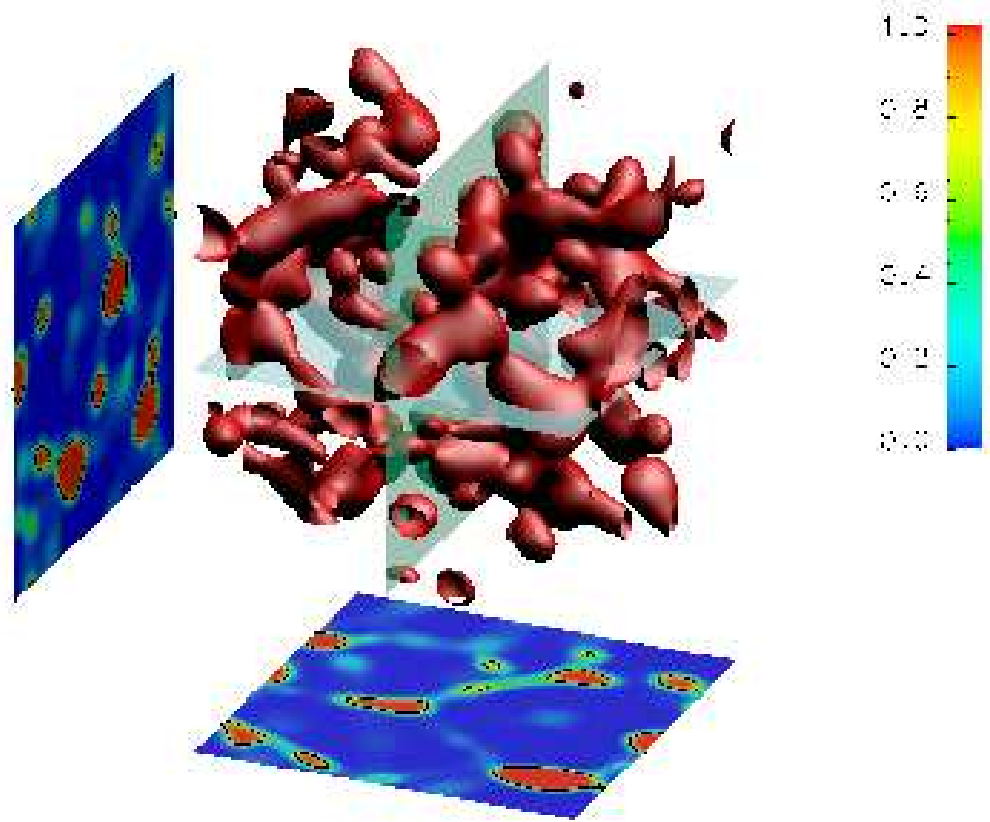}
\caption{Stationary state of the evolution of the states depicted in Figs.
  \ref{fig:3d:in} and \ref{fig:3d:med} according to model
  \eqref{Naldi:modello}. Well developed threedimensional network-like
  structures are observed} 
\miolabel{fig:3d:finale}
\end{figure}

We performed numerical simulations with varying initial average cell
density $\bar{n}$. We observed that the initially randomly
distributed cells coalesce forming elongated structures and evolve
towards a stationary state mimicking the geometry of a blood
vessel network in the early stages of formation.

We assigned $\bar{n}$ in the range $2100-3500\,\mathrm{cells/mm^3}$ and
performed 10 to 15 runs for each density value with a $128\times 128 \times 128$ grid on a biological system of $1~\mathrm{mm}^3$. The characteristic lengths and geometric
properties of  
the stationary state
depend on $\bar{n}$ and we observed a percolative phase transition similar
to the one described in \cite{GAC+03} for the twodimensional case.

\subsection{Analysis of the percolative phase transition}

In experimental blood vessel
formation it has been shown that a percolative transition is observed, by varying the initial
cell density. For low cell densities only isolated clusters of
endothelial cells are observed, while for very high densities cells
fill the whole available space. In between these two extreme
behaviours, close to a critical cell density $n_{c}$,
one observes the formation of critical percolating clusters
connecting opposite sides of the domain,
characterized by well defined scaling laws and exponents. These
exponents are known not to depend on the microscopic details of the
process while their values characterize different classes of aggregation
dynamics.

The purely geometric problem of percolation is actually one of the
simplest phase transitions occurring in nature. Many percolative models
show a second order phase transition in correspondence to a critical
value $n_{\mathrm{c}}$, \textit{i.e.} the probability $\Pi$ of observing
an infinite, percolating cluster is $0$ for $\bar{n}<n_{\mathrm{c}}$ and $1$
for $\bar{n}>n_{\mathrm{c}}$~\cite{SA94}.
The phase transition can be studied by focusing on the values of an
order parameter, i.e. an observable quantity that is zero before the
transition and takes on values of order $1$ after it. In a percolation
problem the natural order parameter is the probability $P$ that
a randomly chosen site belongs to the infinite cluster (on finite
grids, the infinite cluster is substituted by the largest one).

In the vicinity of the critical density $n_{c}$ the geometric
properties of clusters show a peculiar scaling behavior.  For
instance, in a system of linear finite size $L$, the probability of
percolation $\Pi(n,L)$, defined empirically as the fraction of
computational experiments that produce a percolating cluster, is
actually a function of the combination $(n-n_{\mathrm{{c}}})\,
L^{1/\nu}$, where $\nu$ is a universal exponent. 

In a neighborhood of the critical point and on a system of finite size
$L$, the following finite size scaling relations are also observed:

\begin{equation}
\miolabel{pifss}
\Pi(\bar{n},L) \>\sim\>  \widehat{\Pi}[(\bar{n}-n_\mathrm{c})L^{1/\nu}]
\end{equation}

There are two main reasons to study percolation in
relation to vascular network formation: (\emph{i}) percolation is a
fundamental property for vascular networks: blood should have the
possibility to travel through the whole vascular network to carry
nutrients to tissues; (\emph{ii}) critical exponents are robust
observables characterizing the aggregation dynamics.

A rather complete characterization of percolative exponents in the
two-dimensional case has been provided in \cite{GAC+03}.

As a first step in the study of the more realistic threedimensional
case, we compute the exponent $\nu$ characterizing the structures
produced by the model dynamics \eqref{Naldi:modello} with varying
initial cell density.

To this aim, extensive
numerical simulation of system (\ref{Naldi:modello}) were performed
using lattice sizes $L=1,0.78,0.62,0.5\,\mathrm{mm}$, with different
values of the initial density $\bar{n}$.
For each point 10 to 15 realizations of the system of size
$1\mathrm{mm}$ were computed, depending on the proximity to the
critical point.

The continuous density at final time $n(\mathbf{x})$ was then mapped to a set of
occupied and empty sites by choosing a threshold $n_0$. Each region of 
adjacent occupied sites (cluster) was marked with a different
index. The percolation probability $\Pi$ for each set
of realizations was then measured. In Fig.~\ref{fig:clusters} we show
clusters obtained in a box with $L=0.5\,\mathrm{mm}$ with
$\bar{n}=3100$. 
The largest percolating cluster is shown in red,
together with some other smaller clusters shown in different colors.

Using relation \eqref{pifss}, we estimate the position of the critical
point $n_{\mathrm{c}}$ and the value of the critical exponent
$\nu$. The data for different box side
length and initial density should lie on a single curve after rescaling the densities as
$\hat{n}=(\bar{n}-n_c)L^{1/\nu}$. For
fixed $n_{\mathrm{c}}$ and $\nu$ we rescale $\bar{n}$ and
fit the data with a logistic curve, then compute the distance of the
data from the curve. The squared distance is minimized to obtain estimates for $n_{\mathrm{c}}$ and $\nu$.

Using $n_0=0.35$ we obtain $n_{\mathrm{c}}=2658$ and
$\nu=0.84$. This latter value is compatible with the known value $0.88$
for random percolation in three dimensions \cite{SA94}.

\begin{figure}
\includegraphics[width=\textwidth]{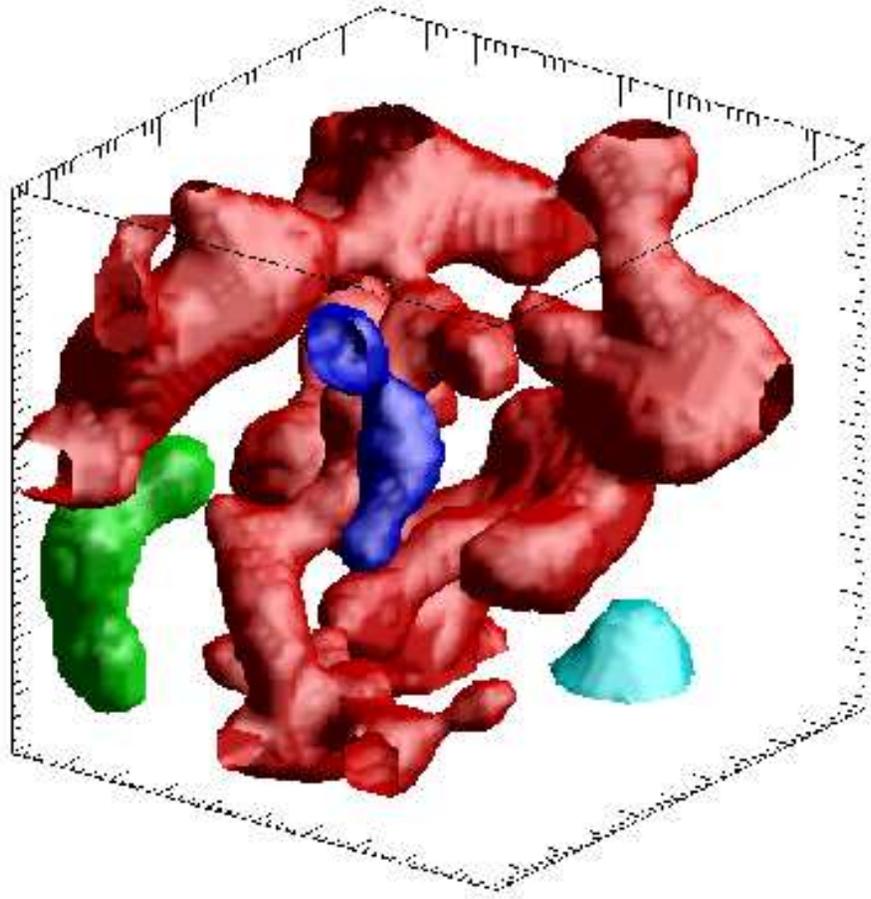}
\caption{Cluster percolation with cell density $n=2500$cells/mm$^{3}$. \textbf{A}:
connected clusters in a realization of model \eqref{Naldi:modello}. \textbf{B}: the largest
cluster depicted in A percolates.}
\miolabel{fig:clusters}
\end{figure}

\begin{figure}
\begin{center}
\includegraphics[width=.95\textwidth]{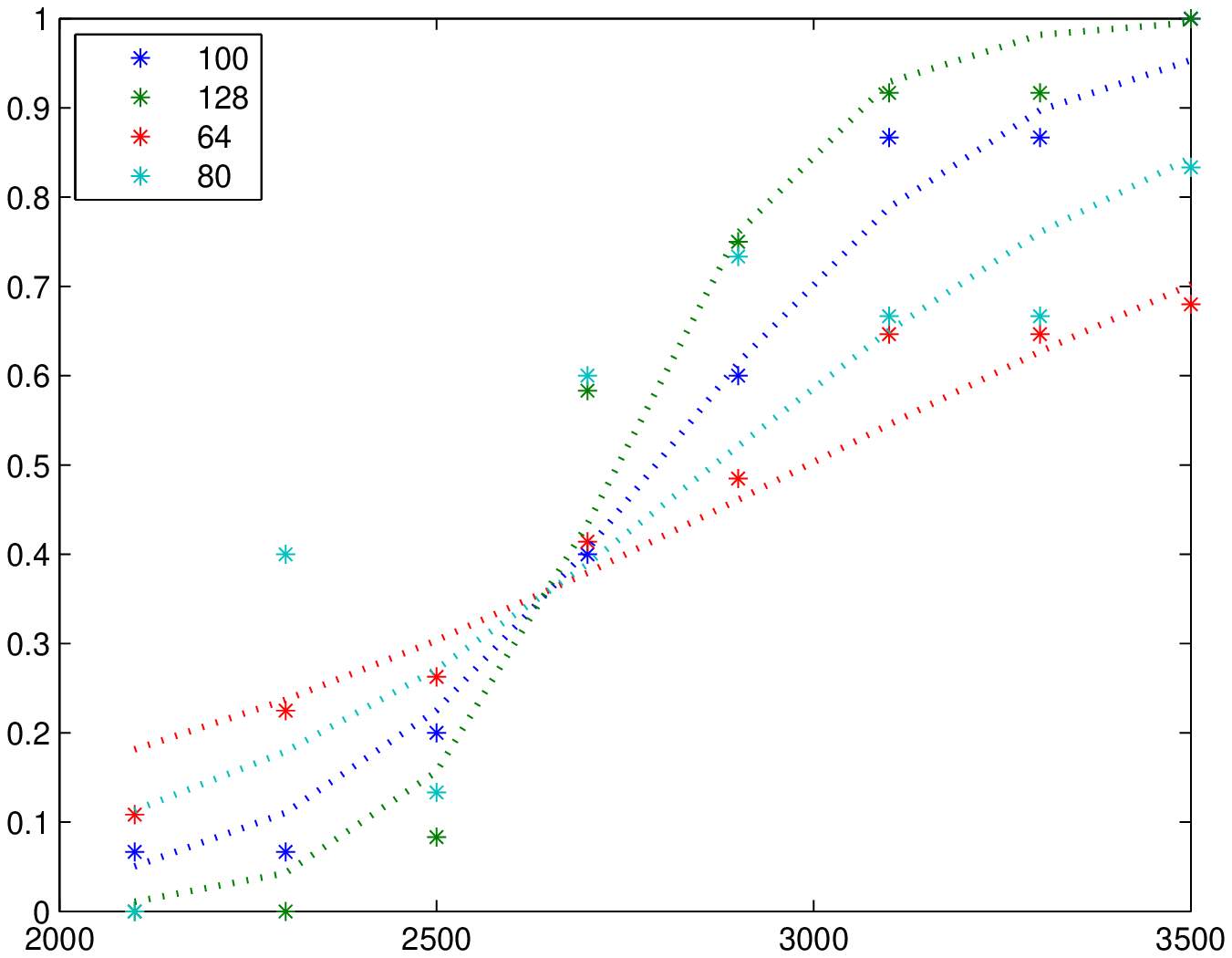}
\includegraphics[width=.95\textwidth]{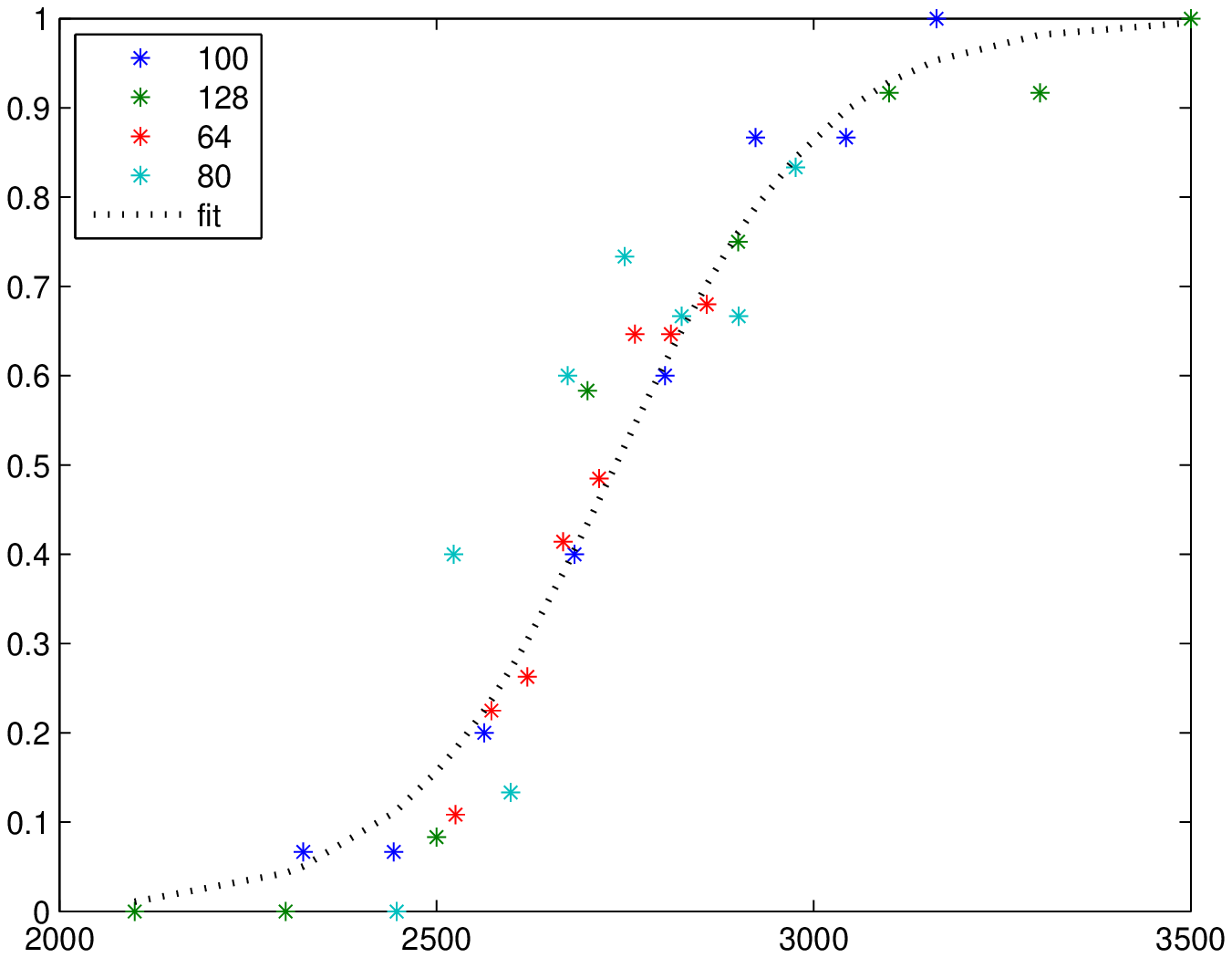}
\end{center}
\caption{Percolation probability at varying densities A: B: the curves in
A are collapsed according to formula \eqref{pifss}}
\miolabel{fig:perc:collapse}
\end{figure}

\section{Conclusions}

We have exposed results on the numerical simulation of vascular
network formation in a threedimensional setting.

We have used the threedimensional version of
the equations proposed in \cite{GAC+03,SAG+03} as a computational model. Evolution starting
from initial conditions mimicking the experimentally observed ones
produce network-like structures qualitatively similar to those
observed in the early stages of in vivo vasculogenesis.

As a starting point towards a quantitative comparison between
experimental data and the theoretical model we nedd to select a set of
observable quantitaties which provide robust quantitative information
on the network geometry.
The lesson learned from the study of twodimensional vasculogenesis
is that percolative exponents are an interesting set of such observables,
so we tested the computation of percolative exponents on simulated
network structures.

A quantitative comparison of the geometrical properties of experimental
and computational network structures will become possible as soon as
an adequate amount of experimental data, allowing proper statistical
computation, will become available.

In order to compute the robust statistical observables described in
the paper one has to perform many runs of the simulation code using
different random initial data. This, toghether with the intensive use
of computational resources required by a three-dimensional
hydrodynamic simulation on fine grids, calls for an efficient
implementation of the computational model on parallel computers, as
the one we presented in this paper.

Simulations of blood vessel structures can in principle present practical
implications. Normal tissue function depends on adequate supply of
oxygen through blood vessels. Understanding the mechanisms of
formation of blood vessels has become a principal objective of medical
research, because it would offer the possibility of testing medical
treatments \textit{in silicio}.  One can think that the dynamical
model (\ref{Naldi:modello}) can be also exploited  in the future to design
properly vascularized artificial tissues by controlling the
vascularization process through appropriate signaling
patterns.

\bibliographystyle{plain}
\bibliography{0604606_v2}

\end{document}